\documentclass[12pt]{article}

\usepackage{xcolor,graphicx}
\usepackage{amssymb,amsmath}

\usepackage{enumerate}

\DeclareMathAlphabet{\eufrak}{U}{}{}{}
\SetMathAlphabet\eufrak{normal}{U}{euf}{m}{n}
\SetMathAlphabet\eufrak{bold}{U}{euf}{b}{n}

\numberwithin{equation}{section}

\newenvironment{Proof}{\removelastskip\par\medskip
\noindent{\em Proof.} \rm}{\penalty-20\null\hfill$\square$\par\medbreak}

 \def\real{{\mathord{\mathbb R}}}
 
 \def\inte{{\mathord{\mathbb N}}}
 
 \def\qu{{\mathord{\mathbb Z}}}

 \def\Dom{{\mathrm{{\rm Dom}}}}

 \def\real{{\mathord{{\rm I\kern-3pt R}}}}        
 
 \def\inte{{\mathord{{\rm I\kern-3pt N}}}}
 \def\sZZ{{\rm Z\kern-.45em{}Z}}
 \def\sQQ{{\kern 0.27em \vrule height1.45ex width0.03em depth0em
           \kern-0.30em \rm Q}}
 \def\qu{{\mathchoice
         {\sQQ}
         {\sQQ}
   {\kern 0.225em \vrule height1.05ex width0.025em depth0em \kern-0.25em \rm Q}
   {\kern 0.180em \vrule height0.78ex width0.020em depth0em \kern-0.20em \rm Q}
         }}
 \def\sGG{{\kern 0.27em \vrule height1.45ex width0.03em depth0em
           \kern-0.30em \rm G}}
 \def\gg{{\mathchoice
         {\sGG}
         {\sGG}
   {\kern 0.225em \vrule height1.05ex width0.025em depth0em \kern-0.25em \rm G}
   {\kern 0.180em \vrule height0.78ex width0.020em depth0em \kern-0.20em \rm G}
         }}

 \newtheorem{prop}{Proposition}[section]
 \newtheorem{lemma}[prop]{Lemma}

 \newtheorem{remark}[prop]{Remark}

 \def\Dom{{\mathrm{{\rm Dom \! \ }}}}

\def\E{\mathop{\hbox{\rm I\kern-0.20em E}}\nolimits}

\def\e{\mathop{\hbox{\rm e}}\nolimits}

\usepackage{empheq} 

\usepackage[colorlinks=true, urlcolor=blue,linkcolor=blue, citecolor=blue]{hyperref}

 \newcounter{hyp}
\title{\huge Selling at the ultimate maximum in a regime-switching model }

\author{
 Yue Liu \qquad
 Nicolas Privault
\\
\vspace{-0.1cm}
\small
School of Physical and Mathematical Sciences
\\
\vspace{-0.1cm}
\small
Division of Mathematical Sciences
\\
\vspace{-0.1cm}
\small
 Nanyang Technological University \\
\vspace{-0.1cm}
\small
Singapore 637371
}

\allowdisplaybreaks

\usepackage{subcaption}

\begin{document}

\hyphenation{func-tio-nals}
\hyphenation{Privault}

\maketitle

\vspace{-1cm}

\begin{abstract}
{ 
 This paper deals with optimal prediction in a regime-switching model driven by a continuous-time Markov chain. We extend existing results for geometric Brownian motion by deriving optimal stopping strategies that depend on the current regime state, and prove a number of continuity properties relating to optimal value and boundary functions. Our approach replaces the use of closed form expressions, which are not available in our setting, with PDE arguments that also simplify the approach of \cite{dutoit} in the classical Brownian case. 
}
\end{abstract}


\baselineskip0.7cm

\section{Introduction}
 Regime-switching models have been introduced by Hamilton~\cite{hamilton}
 in discrete time
 and are among the most popular and effective risky asset models.
 The regime-switching property is reflected in the changes of states 
 of a Markov chain $\beta_t$, which stands for the influence of external 
 market factors.
\\ 
 
 European options have been priced in continuous-time
 regime-switching models by Yao, Zhang and Zhou~\cite{yaozhangzhou} 
 via a recursive algorithm, and in 
 Liu, Zhang and Yin~\cite{liu-zhang-yin} using the 
 fast Fourier transform. 
 Optimal stopping for option pricing 
 in regime-switching models has been considered in 
 Guo~\cite{guoxin}, 
 Guo and Zhang~\cite{guozhang}, 
 Le and Wang~\cite{le-wang}, 
 and in Ly Vath and Pham~\cite{vathana} with optimal switching. 
 Optimal selling under threshold rules has been dealt with 
 in Eloe et al.~\cite{eloe} in an exponential Gaussian diffusion model 
 with regime switching. 
 We refer to Shiryaev~\cite{Shiryaevstopping} and 
 Peskir and Shiryaev~\cite{peskir} 
 for background on the characterization of optimal stopping times
 and rewards.
\\ 
 
 The problem of selling a stock at the ultimate maximum has been considered 
 by Du Toit and Peskir~\cite{dutoit} 
 as an extension of the results of Shiryaev, Xu and 
 Zhou~\cite{shiryaev-xu-zhou}.
 In this paper we extend the result of \cite{dutoit} 
 to the framework of Markovian regime switching. 
 Some of our results are natural extensions of those of \cite{dutoit} 
 by averaging over the regime-switching component, 
 however the regime-switching case presents notable 
 differences and additional difficulties compared with the 
 classical Brownian case. 
 For example, the optimal boundary functions 
 depend on the regime state of the system, and they may not be monotone 
 if the drift coefficients have switching signs.  
 In addition
 we can no longer rely on closed form expressions as in \cite{dutoit} 
 and instead we use PDE arguments, 
 cf. e.g. Lemma~\ref{equationg}, 
 that also simplify the 
 original approach. 
\\ 
 
 In Lemma~\ref{1.1} we write the optimal value of the problem as 
 a function of time, the regime state, and 
 the relative maximum of the underlying asset. 
 In the general case of real-valued drifts
 $\mu (i) \in \real$, $i \in {\cal M}$, 
 we identify the optimal stopping
 time $\tau_D$ 
 in Proposition~\ref{3.1}, and in 
 Proposition~\ref{3d} we determine the structure of the 
 optimal stopping set via 
 its boundary functions $b_{D}(t,j)$ for $i$ in the state 
 space ${\cal M}$ of the regime-switching chain. 
\\ 
 
When the drift parameters $(\mu (i))_{i \in {\cal M}}$
of the regime-switching chain 
 are nonnegative we prove the continuity 
 and monotonicity of boundary functions $b_{D}(t,j)$ 
 in Proposition~\ref{boundarydec}, 
 by extending arguments of \cite{dutoit} 
 to the regime-switching setting.
 Those results are illustrated
 in Figures~\ref{fig2-3} and \ref{fig3} by the plotting of
 value functions that yield the optimal
 stopping boundaries. 
 \\
 
 In Proposition~\ref{2spcase} we show that immediate exercise 
 is optimal when all drift parameters 
 $\mu(i)$ are negative, $i\in{\cal M}$, 
 while exercise at maturity becomes optimal 
 when $\mu(i)\geq \sigma^2(i)$ for all $i\in{\cal M}$, 
 where $\sigma (i)$ are the volatility parameters. 
\\
  
 In Proposition~\ref{bdry} we derive 
 a Volterra type integral equation \eqref{ebdry} which is
 satisfied by the boundary function $b_{D}(t,j)$ of the 
 stopping set. 
 Such an equation is difficult to solve because,
 unlike in the classical setting \cite{dutoit},
 it also relies on the knowledge of the 
 optimal value function, cf. Remark~\ref{r1}. 
 In addition the associated free boundary problem \eqref{lv}-\eqref{lv2}
 consists in a system of interacting PDEs that cannot be
 solved without additional assumptions,
 cf. e.g. Buffington and Elliott~\cite{buffington}
 for a solution under an ordering
 condition on boundary functions in the case of American options. 
 \\ 
 
 A treatment of drifts coefficients $(\mu (i))_{i \in{\cal M}}$ 
 with switching signs has been proposed in 
 of \cite{liu-privault2} via a recursive algorithm
 that does not rely on a Volterra equation.
 In this case it turns out that 
 the boundary functions $b_{D}(t,i)$ may not be decreasing 
 in $t \in [0,T]$. 
 \\

 We proceed as follows. 
 In Section~\ref{s1} we formulate the optimal prediction problem 
 using optimal value functions. 
 In Section~\ref{ss2} we derive the optimal stopping strategies 
 in terms of the hitting time of the boundary function 
 of a stopping set. 
 Section~\ref{ss3} is devoted to continuity lemmas,
 which are used to prove the continuity of boundary functions. 
 In Section~\ref{ss4} we also derive 
 the Volterra integral equation which is satisfied by
 the boundary functions when the drift coefficients are nonnegative. 
 Finally we study particular exercise strategies and we 
 present a numerical simulation of boundary functions. 
\section{Problem formulation} 
\label{s1} 
 Given a standard Brownian motion $(B_t)_{t\in \real_+}$ independent 
 of the Markov chain $(\beta_t)_{t\in \real_+}$,
 and the filtration $({\cal F}_t)_{t\in \real_+}$ generated by $(B_t)_{t\in \real_+}$ and $(\beta_t)_{t\in \real_+}$, we consider
an asset price $(Y_t)_{t\in \real_+}$ modeled
 by a geometric Brownian motion
\begin{equation}
\label{1gb}
dY_t=\mu(\beta_t)Y_tdt+\sigma (\beta_t)Y_tdB_t, \qquad 0 \le t\le T,
\end{equation}
 with regime switching
 driven by a finite-state continuous-time Markov chain
 $(\beta_t)_{t\in \real_+}$ on ${\cal M}=\{1,2, \ldots ,m\}$,
 where 
 $\mu : {\cal M} \longrightarrow \real$, and
 $\sigma : {\cal M} \longrightarrow (0,\infty )$
 are deterministic functions.
 In this paper we deal with the optimal prediction problem
\begin{equation}
\label{vt}
V_t=\inf_{t\le\tau\le T}E\left[\sup_{0\le s\le T }\frac{Y_s}{Y_{\tau}}\Big|{\cal F}_t\right],
\end{equation}
introduced in \cite{dutoit} for geometric Brownian motion,
in which the infimum of expected values 
over all $({\cal G}_s^t)_{s\in [t,T]}$-stopping times $\tau$
  minimizes the ``regret'' of the stopping decision,
  where 
  ${\cal G}_s^t : = \sigma ( B_r - B_t, \ \beta_r \ : \ t \leq r \leq s )$, 
  $s\in [t,T]$. 
\\ 
 
 The next Lemma~\ref{1.1} shows that the optimal value function $V_t$ in \eqref{vt} can be written as a function of $\big( t, \beta_t, \hat{Y}_{0,t} / Y_t \big)$, where $\hat{Y}_{0,t}$ is defined by 
\begin{equation}\label{yts}
\hat{Y}_{t,s} := \max_{t \leq r \leq s} Y_r,
\qquad
0\leq t \leq s \leq T. 
\end{equation}
\vspace{-1cm} 
\begin{lemma}
\label{1.1}
 The optimal value function $V_t$ in \eqref{vt}
 takes the form
\begin{equation*}\label{lem1.1}
V_t =V \big( t, \hat{Y}_{0,t} / Y_t ,\beta_t\big),
\end{equation*}
 where the function 
 $V : [0,T]\times [1,+\infty)\times {\cal M}\rightarrow \real_+$
 is given by
\begin{equation}
\label{vtjx}
V(t,x,j)=\inf_{t\le\tau\le T}E\left[\frac{1}{{Y_\tau}}\max\left( 
 xY_t,\hat{Y}_{t,T}\right)\;\Big|\;\beta_t=j\right],
\end{equation}
 $0 \leq t \leq T$, $x\geq 1$, $j\in {\cal M}$.
\end{lemma}
\begin{Proof}
 Given $t\in[0,T]$, using the drifted Brownian motion
 \begin{equation*}
   \hat {B}^t_u:=B_{u+t}-B_t+\int_{t}^{t+u} \left(\frac{\mu(\beta_s)}{\sigma(\beta_s)}-\frac{\sigma(\beta_s)}{2}\right)ds,
 \quad u\in[0,T-t],
\end{equation*}
 we rewrite the solution of \eqref{1gb} as
\begin{equation}\label{s}Y_s=Y_t\exp \left(
 \int_0^{s-t}\sigma(\beta_{u+t})d\hat {B}^t_u \right),
 \qquad s \in [ t , T],
\end{equation}
 and define
 \begin{equation*}
   \hat {S}^t_s:=\sup_{0\le r\le s}\int_0^r\sigma(\beta_{u+t})d\hat{B}^{t}_{u},
 \qquad s \in [0,T-t].
\end{equation*} 
 By the definition of $\hat Y_{t,s}$ in \eqref{yts} and expression \eqref{s},
 and from the conditional independence of $\left(\big(\hat{B}_{r-t}^t\big)_{r\geq t},
 \big(\hat{S}_{r-t}^t\big)_{r\geq t}\right)$ with ${\cal F}_t$ given $\beta_t$ 
 we have, for any $({\cal G}_s^t)_{s\in [t,T]}$-stopping
 time $\tau$ with values in $[t,T]$, 
 letting $a\vee b = \max(a,b)$, 
 \begin{eqnarray*}
\lefteqn{E\left[\sup_{0\le s\le T }\frac{Y_s}{Y_{\tau}}\Big|{\cal F}_t\right]=E\left[\frac{\hat {Y}_{0,t}}{Y_\tau}
 \vee \frac{\hat{Y}_{t,T}}{Y_\tau}\Big|{\cal F}_t\right]}\nonumber\\
&=&E\left[\left(\frac{\hat {Y}_{0,t}}{Y_t}\; \e^{-\int_0^{\tau-t}\sigma(\beta_{u+t})d\hat {B}^t_u}\right)\vee \e^{\hat {S}^t_{T-t}-\int_0^{\tau-t}\sigma(\beta_{u+t})d\hat {B}^t_u}\Big|{\cal F}_t\right]\nonumber\\
&=&E\left[\left(\frac{\hat {Y}_{0,t}}{Y_t}\;\e^{-\int_0^{\tau-t}\sigma(\beta_{u+t})d\hat {B}^t_u}\right)\vee \e^{\hat {S}^t_{T-t}-\int_0^{\tau-t}\sigma(\beta_{u+t})d\hat {B}^t_u}\;\Big|\;\beta_t,\frac{\hat {Y}_{0,t}}{Y_t}\right]\nonumber\\
&=&E\left[\frac{\hat {Y}_{0,t}\vee\hat{Y}_{t,T}}{Y_\tau}\;\Big|\;\beta_t , \frac{\hat{Y}_{0,t}}{Y_t} \right]\nonumber\\
&=&E\left[\frac{1}{{Y_\tau}}\max\left(xY_t,\hat{Y}_{t,T}\right)\;\Big|\;\beta_t \right]_{x=\hat{Y}_{0,t}/Y_t}
\end{eqnarray*}
where the last line follows from the conditional independence between $\hat{Y}_{0,t}/Y_t$ and
$$\left(\frac{Y_\tau}{Y_t},\frac{\hat{Y}_{t,T}}{Y_\tau}\right)=\left(\exp\left(\int_0^{\tau-t}\sigma(\beta_{u+t})d\hat {B}^t_u\right),\exp\left(\hat {S}^t_{T-t}-\int_0^{\tau-t}\sigma(\beta_{u+t})d\hat {B}^t_u\right)\right)$$
given $\beta_t$. Therefore by definition \eqref{vt} and  expression \eqref{vtjx}, we obtain
\begin{eqnarray*}
V_t&=&\inf_{t\le\tau\le T}E\left[\sup_{0\le s\le T }\frac{Y_s}{Y_{\tau}}\Big|{\cal F}_t\right]\nonumber\\
&=&\inf_{t\le\tau\le T}E\left[\frac{1}{{Y_\tau}}\max\left(xY_t,\hat{Y}_{t,T}\right)\;\Big|\;\beta_t \right]_{x=\hat{Y}_{0,t}/Y_t}\\
&=&V \left( t\frac{\hat{Y}_{0,t}}{Y_t} ,\beta_t,\right).
\end{eqnarray*}
\vskip-1cm
\end{Proof}
 In the next lemma we rewrite the optimal stopping problem \eqref{vt} in the standard form \eqref{vtj} below, using the function 
\begin{equation}
\label{g2}
 G(t,x,i) : =
 E\left[\max \left( x , \hat{Y}_{t,T} / Y_t \right) 
 \;\Big|\;\beta_t = i\right], 
 \quad t\in[0,T], \ i\in{\cal M}, \ x\geq 1, 
\end{equation} 
 with $G(T,x,i)=x$, $x\geq 1$. 
\begin{lemma} 
\label{v=g}
 The function
 $V : [0,T]\times [1,+\infty)\times {\cal M}\rightarrow \real_+$
 defined by \eqref{lem1.1} admits the expression
\begin{equation}
\label{vtj}
 V(t,x,j) =
 \inf_{t\le \tau \le T}
 E \left[
 G \left( \tau, X_\tau^{t,x} ,\beta_\tau\right)\; \Big| \; \beta_t = j
 \right],
\end{equation}
for $t\in[0,T], j\in{\cal M}, x\geq 1$, where
\begin{equation} 
\label{x1}
 X_r^{t,x}:=
 \frac{1}{Y_r}
 \max \left( x Y_t , \hat{Y}_{t,r} \right),
 \qquad r \in [t, T],
 \quad x\geq 1.
\end{equation}
\end{lemma}
\begin{Proof} By a conditional independence argument as in the proof of Lemma~\ref{lem1.1}, for any $s\in[t,T]$ we have
\begin{eqnarray}\label{replace}
E\left[\frac{\hat{Y}_{0,T}}{Y_s}\;\Big|\;{\cal F}_s\right]&=&E\left[\frac{\hat{Y}_{0,s}\vee\hat{Y}_{s,T}}{Y_s}\;\Big|\;{\cal F}_s\right]\nonumber\\
&=&E\left[\frac{\hat{Y}_{0,s}\vee\hat{Y}_{s,T}}{Y_s}\;\Big|\;
\frac{\hat{Y}_{0,s}}{Y_s},\beta_s\right]\nonumber\\
&=&E\left[y\vee\frac{\hat{Y}_{s,T}}{Y_s}\;\Big|\;\beta_s\right]_{y=\hat{Y}_{0,s}/Y_s}\nonumber\\
&=&G\left(s,\frac{\hat{Y}_{0,s}}{Y_s},\beta_s\right).
\end{eqnarray}
Next, we extend the above relation \eqref{replace} to $({\cal G}_s^t)_{s\in [t,T]}$-stopping times $\tau$ written as the limit 
of a decreasing sequence of discrete stopping times by checking that for any
discrete $({\cal G}_s^t)_{s\in [t,T]}$-stopping time
$\tau=\sum_{i=1}^n s_i {\bf 1}_{\{\tau=s_i\}}$,
$s_1,...,s_n\in [t,T]$, $n\geq 1$, by \eqref{replace} we have
\begin{eqnarray*}
E\left[\frac{\hat{Y}_{0,T}}{Y_\tau}\;\Big|\;{\cal F}_\tau\right]&=&\sum_{i=1}^nE\left[\frac{\hat{Y}_{0,T}}{Y_\tau}{\bf 1}_{\{\tau=s_i\}} \;\Big|\;{\cal F}_\tau\right]\nonumber\\
&=&\sum_{i=1}^nE\left[\frac{\hat{Y}_{0,T}}{Y_{s_i}}{\bf 1}_{\{\tau=s_i\}} \;\Big|\;{\cal F}_{s_i}\right]\\
&=&\sum_{i=1}^nE\left[\frac{\hat{Y}_{0,T}}{Y_{s_i}}\;\Big|\;{\cal F}_{s_i}\right]{\bf 1}_{\{\tau=s_i\}}\\
&=&\sum_{i=1}^nG\left(s_i,\frac{\hat{Y}_{0,s_i}}{Y_{s_i}},\beta_{s_i}\right){\bf 1}_{\{\tau=s_i\}}\\
&=&G\left(\tau,\frac{\hat{Y}_{0,\tau}}{Y_\tau},\beta_\tau\right).
\end{eqnarray*}
Taking the conditional expectation $E[\;\cdot\;|\;\beta_t=j,{\hat{Y}_{0,t}}/{Y_t}=x]$ on both sides of the above equality, we obtain
\begin{equation}\label{lem2b}E\left[\frac{\hat{Y}_{0,T}}{Y_\tau}\;\Big|\;\beta_t=j,\frac{\hat{Y}_{0,t}}{Y_t}=x\right]
=E\left[G\left(\tau,\frac{\hat{Y}_{0,\tau}}{Y_\tau},\beta_\tau\right)\;\Big|\;\beta_t=j,\frac{\hat{Y}_{0,t}}{Y_t}=x\right].
\end{equation}
By \eqref{x1} and the conditional independence between ${\hat{Y}_{0,t}}/{Y_t}$ and $\left( 
 {Y_t}/{Y_\tau},{\hat{Y}_{t,\tau}}/{Y_\tau} \right)$ 
 given \pretolerance=369 $\beta_t=j$ we find
 \begin{eqnarray}
   \nonumber 
\lefteqn{\!\!\!\!\!\!\!\!\!\!\!\!\!\!\!\!\!\!\!\!\!\!\!E\left[\frac{1}{{Y_\tau}}\max\left(xY_t,\hat{Y}_{t,T}\right)\;\Big|\;\beta_t=j\right]
=E\left[\frac{\hat{Y}_{0,T}}{Y_\tau}\;\Big|\;\beta_t=j,\frac{\hat{Y}_{0,t}}{Y_t}=x\right]}\nonumber\\
&&\!\!\!\!\!\!\!\!\!\!\!\!\!=E\left[G\left(\tau,\frac{\hat{Y}_{0,\tau}}{Y_\tau},\beta_\tau\right)\;\Big|\;\beta_t=j,\frac{\hat{Y}_{0,t}}{Y_t}=x\right]
\nonumber\\
&&\!\!\!\!\!\!\!\!\!\!\!\!\!=E\left[G\left(\tau,\frac{(xY_t)\vee\hat{Y}_{t,\tau}}{Y_\tau},\beta_\tau\right)
\;\Big|\;\beta_t=j,\frac{\hat{Y}_{0,t}}{Y_t}=x\right]\nonumber\\
&&\!\!\!\!\!\!\!\!\!\!\!\!\!=E\left[G\left(\tau,X_\tau^{t,x},\beta_\tau\right)\;\Big|\;\beta_t=j\right],
\end{eqnarray}
which completes the proof by \eqref{vtjx}.
\end{Proof}
\section{Stopping set and boundary functions} 
\label{ss2} 
 In this section we apply Corollary~2.9 in \cite{peskir}
 in the framework of the regime-switching model \eqref{1gb}
 with $\mu (i) \in \real$, $i\in {\cal M}$, 
 in order to specify
 the stopping set and optimal stopping time associated to the optimal stopping problem \eqref{vt},
 cf. Proposition~\ref{3.1} below.
 In order to deal
 with the existence of an optimal stopping time for \eqref{vt}
 rewritten as \eqref{vtj},
 we define the set
\begin{equation}
\label{ddef}
 D:=\left\{(t, x,j )\in[0,T]\times [1,\infty ) \times {\cal M}\ : \ V(t,x,j)=G(t,x,j)\right\}. 
\end{equation} 
 From the relation $V(T,x,j)=G(T,x,j)=x$, $j\in {\cal M}$, 
 $x\geq 1$, we check that 
 $\{T\}\times [1,\infty ) \times {\cal M} \subset D$, 
 which is consistent with the fact that 
 the infimum in \eqref{vt} is over $({\cal G}_s^t)_{s\in [t,T]}$-stopping times 
 $\tau \in [t,T]$. 
\begin{prop}
\label{3.1}
 Let $t\in[0,T]$.
 Given $\beta_t = j\in {\cal M}$ and ${\hat{Y}_{0,t}}/{Y_t} = x\in[1,\infty)$,
 the $({\cal G}_s^t)_{s\in [t,T]}$-stopping time
\begin{equation}
\label{tau}
 \tau_D(t,x,j) : =\inf\left\{r\geq t \ : \; \left( 
 r,\frac{\hat{Y}_{0,r}}{Y_r},\beta_r \right)\in D\right\}
\end{equation}
 is an optimal stopping time for \eqref{vt}, or equivalently for \eqref{vtj},
 provided it is a.s. finite.
\end{prop}
\begin{Proof} By Corollary~2.9 in \cite{peskir}
 the optimal stopping time for problem \eqref{vtj} exists
 and is equal to $\tau_D (t,x,j)$ in \eqref{tau}
 provided we check that for all $t\in [0,T]$ we have:
\begin{enumerate}[a)] 
\item 
 $G(t,x,j)$ is lower semicontinuous with
 respect to $x$,
 as follows directly from the definition \eqref{g2}
 of $G(t,x,j)$.
\item 
 $V(t,x,j)$ is upper semicontinuous with respect to
 $x$, as follows from the continuity Lemma~\ref{contlemma} 
 below. 
\item We have $E \left[\sup_{t\leqslant s\leqslant T}\left| G(s,X_s^{t,x},\beta_s)\right|\right]<\infty$. 
 Indeed, from \eqref{x1} and \eqref{s} we have
\begin{eqnarray}\label{xtx}
X_s^{t,x}&=& \frac{1}{Y_s}\max \left(xY_t,\hat{Y}_{t,s} \right) \nonumber\\
&=&
 Y_t^{-1} \e^{-\int_0^{s-t}\sigma(\beta_{u+t}) d\hat{B}^t_u}
 \max \left(
 xY_t,Y_t\e^{\hat{S}^t_{s-t}}
 \right)
 \nonumber\\
 &=&
 \e^{\max\left(\log x,\hat{S}^t_{s-t}\right)
 -\int_0^{s-t}\sigma(\beta_{u+t}) d\hat{B}^t_u},
 \quad s\in[t,T].
\end{eqnarray}
Hence by \eqref{g2} and the conditional independence between $X_s^{t,x}=\max\left(xY_t/Y_s,\hat{Y}_{t,s}/Y_s\right)$ and ${\hat{Y}_{s,T}}/{Y_s}$ given $\beta_s$, we find that
\begin{eqnarray}
  \nonumber 
  \lefteqn{G(s,X_s^{t,x},\beta_s)=E\left[y\vee\frac{\hat{Y}_{s,T}}{Y_s}\;\Big|\;\beta_s\right]_{y=X_s^{t,x}}
  } 
\\
\nonumber
&=&E\left[X_s^{t,x}\vee\frac{\hat{Y}_{s,T}}{Y_s}\;\Big|\;\beta_s,X_s^{t,x}\right]
\\
\nonumber
&=&E\left[\e^{-\int_0^{s-t}\sigma(\beta_{u+t})d\hat{B}_u^t}\left(\e^{ \max \left( 
 \log x,\sup_{t\le r\le s}\int_0^{r-t}\sigma(\beta_{u+t})d\hat{B}_u^t\right)}
 \vee \e^{\sup_{s\le r\le T}\int_0^{r-t}\sigma(\beta_{u+t})d\hat{B}_u^t}\right)\;\Big|\beta_s,X_s^{t,x}\right]\\
\label{g+}
&=&E\left[ \e^{ \max \left(\log x,\;\hat{S}_{T-t}^t\right)
-\int_0^{s-t}\sigma(\beta_{u+t})d\hat{B}_u^t }\;|\;\beta_s,X_s^{t,x}\right].
\end{eqnarray}
 Letting 
\begin{eqnarray}
\label{dsjklasdas} 
\check{S}_{T-t}^t:&= & 
\inf_{t\leqslant s \leqslant T}\int_0^{s-t}\sigma(\beta_{u+t})d\hat {B}^t_u\\
\nonumber 
&= & \inf_{t\leqslant s \leqslant T}\left[\int_0^{s-t}\sigma(\beta_{u+t})d{B}_{t+u}+\int_0^{s-t}(\mu(\beta_{u+t})-\sigma^2(\beta_{u+t})/2)du\right],
\end{eqnarray} 
 we conclude that
\begin{align}
E &\left[
 \sup_{t\leqslant s \leqslant T}
 \left| G(s,X_s^{t,x},\beta_s)
 \right|
 \right]
=E\left[\sup_{t\leqslant s \leqslant T} E\left[\e^{\max\left( 
 \log x,\;\hat{S}_{T-t}^t\right) - \int_0^{s-t}\sigma(\beta_{u+t})d\hat {B}^t_u}\;\Big|\;
\beta_s, X_s^{t,x} \right]\right]\nonumber
\\
 &\le
 xE\left[\sup_{t\leqslant s \leqslant T} E\left[\e^{\hat{S}_{T-t}^t-\int_0^{s-t}\sigma(\beta_{u+t})d\hat {B}^t_u}\;\Big|\;
\beta_s, X_s^{t,x} \right]\right]\nonumber
\\
 &\le
 xE\left[\sup_{t\leqslant s \leqslant T} E\left[\e^{\hat{S}_{T-t}^t-\inf_{t\leqslant r \leqslant T}\int_0^{r-t}\sigma(\beta_{u+t})d\hat {B}^t_u}\;\Big|\;
\beta_s, X_s^{t,x} \right]\right]\nonumber
\\
 &=
 x E\left[\e^{\hat{S}_{T-t}^t-\inf_{t\leqslant r \leqslant T}\int_0^{r-t}\sigma(\beta_{u+t})d\hat {B}^t_u}\right]\nonumber
\\
 &=
 x E\left[\e^{\hat{S}_{T-t}^t-\check{S}_{T-t}^t}
 \right]\nonumber
\\
&\leq x
\sqrt{E\left[\e^{2\hat{S}_{T-t}^t}\right]
E\left[\e^{-2\check{S}_{T-t}^t}\right]}
\nonumber\\
&\leq x
\e^{\max_{i\in{\cal M}}|\sigma^2(i)-2\mu(i)|(T-t)}
\sqrt{E\left[\e^{2\hat{S}_{T-t}^t}\right]
E\left[\e^{2\hat{S}_{T-t}^t}\right]
 }
\nonumber\\
 &\leq  x E\left[\e^{2\hat{S}_{T-t}^t}\right]\e^{\max_{i\in{\cal M}}|\sigma^2(i)-2\mu(i)|(T-t)}\nonumber\\
 & < \infty.\label{<infty}
\end{align}
\end{enumerate} 
\vskip-1cm
\end{Proof}
 Define
\begin{equation}
\label{discr}
 F(t,x,j) : = V(t,x,j)-G(t,x,j)\leq 0, 
\end{equation}
 which is nonpositive by \eqref{vtj}, 
 $t\in[0,T]$, $j\in {\cal M}$, $x\geq 1$, 
 so that we have
\begin{equation*}
 D=\left\{(t,x,j)\in[0,T]\times [1,\infty) \times {\cal M} \ : \ F(t,x,j)=0\right\},
\end{equation*}
 hence $D$ is closed from the continuity of $(t,x) \longmapsto V(t,x,j)$ and 
 $(t,x) \longmapsto G(t,x,j)$ on $[0,T] \times [1,\infty )$, 
   cf. Lemmas~\ref{contlemma} and Lemmas~\ref{contlemma1} below,
   respectively. 
 The continuation set $C = D^c$ is an open set that can be written as
$$C=\left\{(t,x,j)\in[0,T]\times [1,\infty) \times {\cal M} \ : \ F(t,x,j)<0\right\}.
$$ 
 In the next Proposition~\ref{3d} 
 we characterize the shape of the stopping set $D$ defined in \eqref{ddef}
 in terms of the boundary function $b_D(t,j)$ defined by 
\begin{equation}\label{bd}b_D (t,j): =\inf\{x \in [1,\infty) \ : \
 (t,x,j)\in D\}. 
\end{equation} 
 From the relation $\{T\} \times [1,\infty ) \times {\cal M} \subset D$ 
 we deduce the terminal condition $b_D (T,j)=1$, $j\in {\cal M}$. 
\begin{prop}
\label{3d}
 For any $(t,x,j)\in[0,T]\times [1,\infty) \times {\cal M}$ such that $(t,x,j)\in D$ we have
\begin{equation} 
\label{fdhsjksdfds} 
\{ t \} \times [x,\infty ) \times \{ j \} \subset D.
\end{equation} 
 and 
\begin{equation}\label{dfinal}D=\left\{
 (t,y,j) \in [0,T] \times [1,\infty)\times {\cal M} 
 \ :\ y\ge b_D (t,j)\right\}.
\end{equation}
\end{prop}
\begin{Proof} 
 Let $y : = \sup \{ z \in [x, \infty ) \ : \ 
 \{ t \} \times [x,z] \times \{ j \} \subset D\}$. 
 If $y<\infty $ then we have $(t,y,j) \in D$ by the closedness
 of $D$, and from the monotonicity property of $F(t,x,j)$ stated
 in Lemma~\ref{3.3}, $(t,y,j) \in D$ admits a right neighborhood 
 of the form 
$$
 \{t \} \times [x,x+\eta] \times \{ j \} \subset D
$$
 for some $\eta>0$, which leads to a contradiction. 
 Hence $y=+\infty$ and \eqref{fdhsjksdfds} holds. 
 Relation \eqref{dfinal} follows from the equivalence 
$$ 
(t,x,j)\in D\Longleftrightarrow \{ t \} \times [x,\infty ) \times \{ j \} \subset D \Longleftrightarrow x\geq b_D(t,j)
$$ 
 that follows from \eqref{bd}. 
\end{Proof}
 The following lemma has been used in the proof of Proposition~\ref{3d}.
\begin{lemma}
\label{3.3}
 For any $(t,x,j)\in D$, we have
\begin{equation}\label{F'>0}
\liminf\limits_{\varepsilon\searrow 0}\frac{F(t,x+\varepsilon ,j)-F(t,x,j)}{\varepsilon} \geq 0.
\end{equation}
\end{lemma}
\begin{Proof} We split the proof into two parts. 
\\ 
$(i)$ From \eqref{g+} we have
\begin{align*}G(s,X_s^{t,x},\beta_s)&=E\left[X_s^{t,x}\vee\frac{\hat{Y}_{s,T}}{Y_s}\Big|\beta_s,X_s^{t,x}\right]
=E\left[X_s^{t,x}\vee\frac{\hat{Y}_{s,T}}{Y_s}\Big|{\cal F}_s\right]\\
&=E\left[ \e^{ \max\left(\log x,\;\hat{S}_{T-t}^t\right)
-\int_0^{s-t}\sigma(\beta_{u+t})d\hat{B}_u^t }\;\Big|\;{\cal F}_s\right],\quad s\in[t,T],
\end{align*}
which extends to any $({\cal G}_s^t)_{s\in [t,T]}$-stopping time $\tau \in[t,T]$ as
\begin{equation}\label{g++}
G(\tau,X_\tau^{t,x},\beta_\tau)
=E\left[ \e^{ \max\left(\log x,\;\hat{S}_{T-t}^t\right)
-\int_0^{\tau-t}\sigma(\beta_{u+t})d\hat{B}_u^t }\;\Big|\;{\cal F}_\tau\right],
\end{equation}
as in \eqref{replace}-\eqref{lem2b} above. For all $x\geq 1$ and $\varepsilon> 0$, consider the $({\cal G}_s^t)_{s\in [t,T]}$-stopping time $$\tau_\varepsilon^+:=\tau_D(t,x+\varepsilon ,j)\in[t,T]$$
defined by \eqref{tau}, which solves the optimal stopping problem
\begin{equation*}
 V(t,x+\varepsilon ,j) =
 \inf_{t\le \tau \le T}
 E \left[
 G \left( \tau, X_\tau^{t,x+\varepsilon},\beta_\tau\right)\; \Big| \; \beta_t = j
 \right]
 = E \left[
 G \left( \tau_\varepsilon^+, X_{\tau_\varepsilon^+}^{t,x+\varepsilon},\beta_{\tau_\varepsilon^+}\right)\; \Big| \; \beta_t = j
 \right],
\end{equation*}
 cf. \eqref{vtj}. The following argument relies on the fact that for any $(t,x,j)\in D$ we have
\begin{equation}
\label{djkldd}
\lim\limits_{\varepsilon\to 0}\tau_D(t,x+\varepsilon ,j)=t,
\end{equation}
 as will be shown in part $(ii)$ below. 
 Relations~\eqref{g2}, \eqref{vtj}, \eqref{g++} and \eqref{djkldd} 
 imply 
\begin{eqnarray}
\lefteqn{\liminf\limits_{\varepsilon\searrow 0} 
\frac{V(t,x+\varepsilon ,j)-V(t,x,j)}{\varepsilon}
}
\nonumber\\
&\geq&\liminf\limits_{\varepsilon\searrow 0} \frac{1}{\varepsilon}E\left[G(\tau_\varepsilon^+,X_{\tau_\varepsilon^+}^{t,x+\varepsilon},\beta_{\tau_\varepsilon^+})
-G(\tau_\varepsilon^+,X_{\tau_\varepsilon^+}^{t,x},\beta_{\tau_\varepsilon^+})\;|\;\beta_t=j\right]\nonumber\\
&=&\liminf\limits_{\varepsilon\searrow 0}\frac{1}{\varepsilon}E\left[E\left[
 \e^{\log (x+\varepsilon)\vee\hat{S}_{T-t}^t-\int_0^{\tau_\varepsilon^+-t}\sigma(\beta_{u+t})d\hat B_u^t}
 -\e^{\log x\vee\hat{S}_{T-t}^t-\int_0^{\tau_\varepsilon^+-t}\sigma(\beta_{u+t})d\hat B_u^t}\Big|{\cal F}_{\tau_\varepsilon^+}\right]\Big|\beta_t=j \right]\nonumber\\
 &=&\liminf\limits_{\varepsilon\searrow 0} \frac{1}{\varepsilon}E\left[
 \e^{\log (x+\varepsilon)\vee\hat{S}_{T-t}^t-\int_0^{\tau_\varepsilon^+-t}\sigma(\beta_{u+t})d\hat B_u^t}
 -\e^{\log x\vee\hat{S}_{T-t}^t-\int_0^{\tau_\varepsilon^+-t}\sigma(\beta_{u+t})d\hat B_u^t}\Big|\beta_t=j \right]\nonumber\\
&=&\liminf\limits_{\varepsilon\searrow 0}\frac{1}{\varepsilon}E\left[
 \e^{\log (x+\varepsilon)\vee\hat{S}_{T-t}^t}
 -\e^{\log x\vee\hat{S}_{T-t}^t}\;\Big|\;\beta_t=j \right]\nonumber\\
&=&\frac{\partial}{\partial x}E\left[
 \e^{\max\left(\log x,\hat{S}_{T-t}^t\right)}\;\Big|\;\beta_t=j \right]\nonumber\\
&=&\frac{\partial G}{\partial x}(t,x,j), 
\label{xup}
\end{eqnarray}
 hence we conclude to \eqref{F'>0}. 
 Here we used the dominated convergence theorem with the bound 
\begin{eqnarray*} 
\lefteqn{ 
 \!\!\!\!\!\!\!\!\! \!\!\!\!\!\!\!\!\! 
 \frac{1}{\varepsilon} 
 \Big|
 \e^{\log (x+\varepsilon)\vee\hat{S}_{T-t}^t-\int_0^{\tau_\varepsilon^+ - t }\sigma(\beta_{u+t})d\hat B_u^t}
 -\e^{\log x\vee\hat{S}_{T-t}^t-\int_0^{\tau_\varepsilon^+ - t }\sigma(\beta_{u+t})d\hat B_u^t} 
 \Big| 
} 
\\
\nonumber
 &\leq & \Big|\frac{\e^{\log(x+\varepsilon)}-\e^{\log x}}{\varepsilon}\Big|\e^{-\inf_{0\leq s\leq T-t}\int_0^s\sigma(\beta_{t+u})d\hat{B}_u^t}=\e^{-\check{S}^t_{T-t}},\label{dominate}
\end{eqnarray*} 
 where $\check{S}_{T-t}^t$ is defined in \eqref{dsjklasdas} 
 and the righthand side is integrable as in the derivation of \eqref{<infty}. 
\\ 
\noindent 
$(ii)$ We turn to the proof of \eqref{djkldd}. 
 From the expression \eqref{vtjx} in Lemma~\ref{lem1.1}, we have
\begin{align}\label{vtjx2}
V(t,x,j)&=\inf_{t\le\tau\le T}E\left[\frac{xY_t\vee \hat{Y}_{t,T}}{{Y_\tau}}\;\Big|\;\beta_t=j,{\cal F}_t\right]\nonumber\\
&=\inf_{t \le\tau\le T}
 E \left[\e^{-\int_0^{\tau-t}\sigma(\beta_{u+t})d\hat{B}_u^t}\big(x\vee \e^{\hat{S}_{T-t}^t}\big)\;\Big|\;\beta_t=j,{\cal F}_t\right]. 
\end{align}
 From \eqref{vtjx2} and
$$ 
 X_r^{t,x+\varepsilon}=\e^{-\int_0^{r-t}\sigma(\beta_{u+t})d\hat B_u^t} 
 (x+\varepsilon \vee \e^{\hat{S}_{r-t}^t})
$$
 cf. \eqref{xtx}, we obtain
\begin{align}\label{vtjx3}
V(r,X_r^{t,x+\varepsilon},\beta_r)&=\inf_{r\le\tau\le T}
 E\left[\e^{-\int_0^{\tau-r }\sigma(\beta_{u+r}) d\hat{B}_u^r }\big(y\vee \e^{\hat{S}_{T-r}^r }\big)\;\Big |\; {\cal F}_r \right]_{y=X_r^{t,x+\varepsilon}}\nonumber\\
&=\inf_{r\le\tau\le T}
 E\left[\e^{-\int_0^{\tau-r }\sigma(\beta_{u+r}) d\hat{B}_u^r }\big( X_r^{t,x+\varepsilon}\vee 
 \e^{\hat{S}_{T-r}^r }\big)\;\Big |\; {\cal F}_r \right].
\end{align}
 Next, from 
 the definition \eqref{tau} of $\tau_D(t,x+\varepsilon ,j)$ and \eqref{vtjx3} we have, on the event $\{\beta_t=j\}$,
\begin{eqnarray}
\label{taueps1}
\lefteqn{
 \tau_D(t,x+\varepsilon ,j) =
 \inf \{r\geq t \ : \ (r,X_r^{t,x+\varepsilon},\beta_r)\in D\}
}
\\
\nonumber
&=&\inf\left\{r\geq t \ : \ \! \inf_{r\le\tau\le T}
 \! \! \!
 E\left[\e^{-\int_0^{\tau-r}\sigma(\beta_{u+r})d\hat{B}_u^r}\big(X_r^{t,x+\varepsilon}\vee \e^{\hat{S}_{T-r}^r}\big)\;\Big |\; {\cal F}_r\right]=
E\left[X_r^{t,x+\varepsilon}\vee \e^{\hat{S}_{T-r}^r }\;\Big |\; {\cal F}_r \right]\right\}
\\
\nonumber
&\leq&\inf\left\{r\geq t \ : \ \! \inf_{r\le\tau\le T}
 \! \! \!
 E\left[\e^{-\int_0^{\tau-r}\sigma(\beta_{u+r})d\hat{B}_u^r} \big(X_r^{t,x}\vee \e^{\hat{S}_{T-r}^r}\big)\;\Big |\;{\cal F}_r \right]\geq 
E\left[X_r^{t,x+\varepsilon}\vee \e^{\hat{S}_{T-r}^r}\;\Big |\;{\cal F}_r \right]\right\}
\\
\nonumber
&\leq&\inf\left\{r\geq t \ : \ \! \inf_{r\le\tau\le T}
 \! \! \!
 E\left[\e^{-\int_0^{\tau-r}\sigma(\beta_{u+r})d\hat{B}_u^r}\big(X_r^{t,x}\vee \e^{\hat{S}_{T-r}^r}\big)\;\Big |\;{\cal F}_r \right]\geq 
\e^{\varepsilon}E\left[X_r^{t,x}\vee \e^{\hat{S}_{T-r}^r} \;\Big |\; {\cal F}_r\right]\right\},
\end{eqnarray}
 where we applied the inequality
$$ 
X_r^{t,x+\varepsilon} 
= \e^{-\int_0^{r-t}\sigma(\beta_{u+t})d\hat{B}_u^t} 
\big(\e^{\log(x+\varepsilon)} \vee \e^{\hat{S}_{r-t}^t}\big)
\leq \e^{-\int_0^{r-t}\sigma(\beta_{u+t})d\hat{B}_u^t+\varepsilon} 
 \big(\e^{\log(x)} \vee \e^{\hat{S}_{r-t}^t}\big)
= \e^{\varepsilon}X_r^{t,x},
$$ 
 $x\geq 1$, $\varepsilon\geq 0$, $r\in[t,T]$. 
 This implies
\begin{eqnarray}
\label{taueps2}
\lefteqn{ 
\lim\limits_{\varepsilon\to 0}\tau_D(t,x+\varepsilon ,j)
} 
\\ 
\nonumber
 & \leq & \lim\limits_{\varepsilon\to 0}\inf\left\{r\geq t \ : \ \! \inf_{r\le\tau\le T}
 E\left[\e^{-\int_0^{\tau-r}\sigma(\beta_{u+r})d\hat{B}_u^r} 
 \big(X_r^{t,x} \vee \e^{\hat{S}_{T-r}^r}\big)\;\Big |\; {\cal F}_r\right]\geq 
 \e^{\varepsilon}E\left[ X_r^{t,x} \vee \e^{\hat{S}_{T-r}^r} \;\Big |\; {\cal F}_r \right]\right\}
\\
\nonumber
&=& \inf\left\{r\geq t \ : \ \! \inf_{r\le\tau\le T}
 \! \! \! E\left[\e^{-\int_0^{\tau-r}\sigma(\beta_{u+r})d\hat{B}_u^r} 
 \big(X_r^{t,x} \vee \e^{\hat{S}_{T-r}^r}\big)\;\Big |\;{\cal F}_r\right]\geq E\left[ X_r^{t,x} \vee \e^{\hat{S}_{T-r}^r} \;\Big |\;{\cal F}_r \right]\right\}
\\
\nonumber
&=& \inf\left\{r\geq t \ : \ \! \inf_{r\le\tau\le T}
 \! \! \!
 E\left[\e^{-\int_0^{\tau-r}\sigma(\beta_{u+r})d\hat{B}_u^r}\big(X_r^{t,x} \vee \e^{\hat{S}_{T-r}^r}\big)\;\Big |\; {\cal F}_r \right]=E\left[ 
 X_r^{t,x} \vee \e^{\hat{S}_{T-r}^r} \;\Big |\; {\cal F}_r\right]\right\}
\\
\nonumber
& = & \inf \left\{r\geq t \ : \ \! (r,X_r^{t,x},\beta_r)\in D\right\}\nonumber\\
& = & t,
\end{eqnarray} 
 since $(t,x,j)\in D$, $\beta_t=j$ and $X_t^{t,x}=x$. Since $\tau_D(t,x+\varepsilon ,j)\geq t$
 we conclude to \eqref{djkldd}.
\end{Proof}
\section{Continuity lemmas} 
\label{ss3} 
 The following property of 
 smooth fit, namely the continuity of the function 
 $y \longmapsto \displaystyle \frac{\partial V}{\partial y} ( t,y,j)$ over the optimal 
 stopping boundary $\partial C$, 
 will be needed in the proof of Proposition~\ref{bdry} below. 
\begin{lemma}\label{sf} 
 For any $(t,y,j)\in \partial C$, $y>1$, we have
$$
 \frac{\partial V}{\partial y} (t,y+,j) 
 =\frac{\partial V}{\partial y} (t,y-,j). 
$$
\end{lemma}
\begin{Proof}
 For any $\varepsilon\in(0,y-1)$, let $\tau^-_\varepsilon=\tau_D(t,y-\varepsilon ,j)\in[t,T]$, cf. \eqref{tau}. Since $(t,y,j)\in \partial C$ and $D$ is closed we have $(t,y,j)\in D$. Similarly to \eqref{taueps1} to \eqref{taueps2}, $\tau^-_\varepsilon$ 
 converges to $t$ a.s. when $\varepsilon$ tends to 0. 
 By the same approach as in \eqref{xup}, 
 replacing $y+\varepsilon$ with $y-\varepsilon$ shows that 
\begin{equation*}
\frac{\partial G}{\partial y}(t,y,j) 
 \leq 
 \liminf_{\varepsilon\searrow 0} \frac{V(t,y-\varepsilon ,j)-V(t,y,j)}{\varepsilon}.
\end{equation*}
 On the other hand, since $(t,y,j)\in \partial C\subset D$, we have 
$$ 
 \limsup_{\varepsilon\searrow 0} \frac{V(t,y-\varepsilon ,j)-V(t,y,j)}{\varepsilon}
 \leq 
 \lim\limits_{\varepsilon\searrow 0} \frac{G(t,y-\varepsilon ,j)-G(t,y,j)}{\varepsilon} 
 = \frac{\partial G}{\partial y} (t,y,j), 
$$
 hence
$$ 
 \frac{\partial V}{\partial y}(t,y-,j) = \frac{\partial G}{\partial y} (t,y,j).$$
 Finally the fact that $V=G$ on the closed set $D$ implies 
$$\frac{\partial V}{\partial y}(t,y-,j)=\frac{\partial V}{\partial y}(t,y+,j)=\frac{\partial G}{\partial y}(t,y,j). 
$$
\end{Proof}
 In the next proposition, 
 which will be used in the proof of Proposition~\ref{bdry}, 
 we show the normal reflection of the free boundary problem
 by proving that
 the right derivative of the value function $V(t,y,j)$ vanishes at $y=1$, 
 cf. also page~264 of \cite{peskir} without regime switching. 
\begin{lemma} 
\label{nr}
 For any $t\in[0,T]$ and $j\in {\cal M}$ we have
$$
 \frac{\partial V}{\partial y}(t,1+,j)=0.
$$
\end{lemma}
\begin{Proof}
 For convenience of notation
 we set $\tau_0 =\tau_D(t,1,j)$, and note that
\begin{align*}
& \limsup_{\varepsilon\searrow 0}\frac{V(t,1+\varepsilon ,j)-V(t,1,j)}{\varepsilon}
\\
&\le  \limsup_{\varepsilon\searrow 0}\frac{1}{\varepsilon}E[G(\tau_0 ,X_{\tau_0 }^{t,1+\varepsilon},\beta_{\tau_0 })
-G(\tau_0 ,X_{\tau_0 }^{t,1},\beta_{\tau_0 })\;|\;\beta_t=j]\\
&= \limsup_{\varepsilon\searrow 0}\frac{1}{\varepsilon} E\left[
 \e^{\log(1+\varepsilon)\vee\hat{S}_{T-t}^t-\int_0^{\tau_0 -t}\sigma(\beta_{t+r})d\hat{B}_r^t}
 -\e^{\hat{S}_{T-t}^t-\int_0^{\tau_0 -t}\sigma(\beta_{t+r})d\hat{B}_r^t}\;\Big|\;\beta_t=j \right] 
\\
& = \limsup_{\varepsilon\searrow 0}E\left[\frac{1}{\varepsilon} 
 \left( 
 \e^{\log(1+\varepsilon)\vee\hat{S}_{T-t}^t-\int_0^{\tau_0 -t}\sigma(\beta_{t+r})d\hat{B}_r^t}
 - \e^{\hat{S}_{T-t}^t-\int_0^{\tau_0 -t}\sigma(\beta_{t+r})d\hat{B}_r^t} 
 \right) {\bf 1}_{\{\hat S^t_{T-t}<\log(1+\varepsilon)\}}\Big|\;\beta_t=j\right]\\
& = E\left[\limsup_{\varepsilon\searrow 0}
 \frac{
 \e^{\log(1+\varepsilon)\vee\hat{S}_{T-t}^t-\int_0^{\tau_0 -t}\sigma(\beta_{t+r})d\hat{B}_r^t}-
 \e^{\hat{S}_{T-t}^t-\int_0^{\tau_0 -t}\sigma(\beta_{t+r})d\hat{B}_r^t}}{\varepsilon}
 {\bf 1}_{\{\hat S^t_{T-t}<\log(1+\varepsilon)\}}\Big|\;\beta_t=j\right]\\
&=0,
\end{align*}
\noindent
 since 
 $\lim\limits_{\varepsilon\searrow 0}{\bf 1}_{\{\hat S^t_{T-t}< \log(1+\varepsilon)\}}=0$, 
 where we applied the dominated convergence theorem as in the proof of 
 Lemma~\ref{3.3} with the same dominating function as in \eqref{dominate}. 
 Since $V(t,y,j)$ is nondecreasing in $y\in[1,\infty)$, we have
 $$ 
 \liminf_{\varepsilon\searrow 0}\frac{V(t,1+\varepsilon ,j)-V(t,1,j)}{\varepsilon}\geq 0, 
$$
 which shows that 
$$ \frac{\partial V}{\partial y}(t,1+,j)=\lim\limits_{\varepsilon\searrow 0}\frac{V(t,1+\varepsilon ,j)-V(t,1,j)}{\varepsilon}=0.$$
\end{Proof}
 Next we consider the infinitesimal generator 
\begin{equation}\label{gene2}
\mathbb{L}f(s,x,j)=\left(\frac{\partial}{\partial s}+x(\sigma^2(j)-\mu(j))\frac{\partial}{\partial x}+\frac{1}{2}\sigma^2(j)x^2\frac{\partial^2 }{\partial x^2}\right)f(s,x,j)+\sum_{i=1}^mq_{j,i}f(s,x,i),
\end{equation}
 of the Markov process $(s,X_s^{t,x},\beta_s)_{s\in[t,T]}$, 
 where $Q=(q_{ij})_{i,j=1,\ldots ,m}$ is the infinitesimal matrix generator of 
 the Markov process $(\beta_t)_{t\in[0,T]}$, 
 for any sufficiently differentiable function $f$ of $(s,y,j)\in[0,T]\times[1,\infty)\times{\cal M}$, 
 cf. Lemma~\ref{mgenerator} below. 
\\ 
 
 The following lemmas will 
 be used in the proof of Proposition~\ref{2spcase} below. 
 In Lemma~\ref{equationg} we replace the use of closed form 
 expressions for $\mathbb{L}G(t,x,j)$, which are no longer available 
 in our setting, with the differential expression \eqref{sakjldsadassa}. 
\begin{lemma} 
\label{equationg} 
 We have 
\begin{equation} 
\label{sakjldsadassa1} 
 \displaystyle \frac{\partial G}{\partial x}(t,1+,j) = 0, 
 \qquad 
 t\in [0,T], 
\end{equation} 
 and 
\begin{equation}
\label{sakjldsadassa} 
 \mathbb{L}G(t,x,j) = 
x\sigma^2(j)\frac{\partial G}{\partial x}(t,x,j)-\mu(j)G(t,x,j), 
 \qquad t\in [0,T], 
\end{equation} 
 with $\mathbb{L}G(T,x,j) = -\mu (j) x$, 
 $j\in {\cal M}$, $x\in [1,\infty)$. 
 In particular, 
 for any $(t,x,j)\in [0,T)\times [1,\infty)\times {\cal M}$ 
 we have 
\begin{equation} 
\label{djkladas} 
\left\{
\begin{array}{ll}
\mathbb{L}G(t,x,j) > 0, & \text{when}\; \mu(j)\leq 0, 
\\ 
\\
\mathbb{L}G(t,x,j) < 0, & \text{when}\; \mu(j)\geq \sigma^2(j). 
\end{array}
\right.
\end{equation} 
 In addition, $\mathbb{L}G(t,x,j)$ is nondecreasing and continuous 
 in $t$ for all $x\geq 1$ when $\mu(j)\geq 0$. 
\end{lemma}
\begin{Proof} 
 For all $j \in {\cal M}$ we let 
 \begin{equation}
   \label{f=g}
f(t,y,z,j) : =yG\left(t,\frac{z}{y} ,j\right) 
 =E\left[\max \left( 
 z,y\frac{\hat{Y}_{t,T}}{Y_t} \right) 
 \;\Big|\;\beta_t=j\right], 
 \quad t \in [0,T], \quad y , z >0. 
\end{equation} 
 By \eqref{1gb} and the It\^o formula we have 
 \begin{align}
   \nonumber
   df\big(t,Y_t,\hat{Y}_{0,t},\beta_t\big)=&\frac{\partial f}{\partial t}\big(t,Y_t,\hat{Y}_{0,t},\beta_t\big)dt+\mu(\beta_t)Y_t\frac{\partial f}{\partial x}\big(t,Y_t,\hat{Y}_{0,t},\beta_t\big)dt\\
   \nonumber
&+\sigma(\beta_t)Y_t\frac{\partial f}{\partial x}\big(t,Y_t,\hat{Y}_{0,t},\beta_t\big)dB_t+\frac{1}{2}\sigma^2(\beta_t)Y^2_t\frac{\partial^2 f}{\partial x^2}\big(t,Y_t,\hat{Y}_{0,t},\beta_t\big)dt\\
   \nonumber
   &+\frac{\partial f}{\partial y}\big(t,Y_t,\hat{Y}_{0,t},\beta_t\big)d\hat{Y}_{0,t} 
 + f\big(t ,Y_t,\hat{Y}_{0,t},\beta_t\big) - f\big(t ,Y_t,\hat{Y}_{0,t},\beta_{t^-}\big)
, 
\end{align}
 and given that 
$$ 
 f\big(t,Y_t,\hat{Y}_{0,t},\beta_t\big)=E\big[\;\hat{Y}_{0,T}\;|\;\beta_t, \ Y_t, \ \hat{Y}_{0,t}\;\big]=E\big[\;\hat{Y}_{0,T}\;|\;{\cal F}_t\;\big], \qquad t\in [0,T], 
 $$
 is a martingale 
 and $\big(\hat{Y}_{0,t}\big)_{t\in[0,T]}$ has finite variation, 
 we find 
\begin{equation}\label{bspde2}
\frac{\partial f}{\partial t}(t,y,z,j)+\mu(j)y\frac{\partial f}{\partial x}(t,y,z,j)
+\frac{1}{2}\sigma^2(j)y^2\frac{\partial^2 f}{\partial x^2}(t,y,z,j)+\sum_{i=1}^{m}q_{j,i}f(t,y,z,i)=0, 
\end{equation}
 and $\displaystyle \frac{\partial f}{\partial y}(t,x,y,j)_{x = y} = 0$. 
 Substituting \eqref{f=g} into \eqref{bspde2} shows that 
\begin{align}\nonumber 
  &y\frac{\partial G}{\partial t}\left(t,\frac{z}{y},j\right)+\mu(j)y\left(G\left(t,\frac{z}{y},j\right)+y\frac{\partial G}{\partial x}\left(t,\frac{z}{y},j\right)\left(-\frac{z}{y^2}\right)\right)\\
  \nonumber
  &+\frac{1}{2}\sigma^2(j)y^2\left(\frac{\partial G}{\partial x}\left(t,\frac{z}{y},j\right)\left(-\frac{z}{y^2}\right)+\frac{z}{y^2}\frac{\partial G}{\partial x}\left(t,\frac{z}{y},j\right)+\frac{z^2}{y^3}\frac{\partial^2 G}{\partial x^2}\left(t,\frac{z}{y},j\right)\right)\\
  \nonumber
&+\sum_{i=1}^{m}q_{j,i}yG\left(t,\frac{z}{y},i\right)=0,
\end{align}
 which shows that the function $G(t,x,j)$ satisfies the PDE 
$$\mu(j)G(t,x,j)+\frac{\partial G}{\partial t}(t,x,j)-\mu(j)x\frac{\partial G}{\partial x}(t,x,j)
+\frac{1}{2}\sigma^2(j)x^2\frac{\partial^2 G}{\partial x^2}(t,x,j)+\sum_{i=1}^mq_{j,i}G(t,x,i)=0, 
$$ 
 and we conclude to \eqref{sakjldsadassa} by \eqref{gene2}. 
 Next we note that 
 \eqref{sakjldsadassa1} follows from 
\begin{equation} 
\label{G'x}
\frac{\partial G}{\partial x}(t,x,j)=P\left(\frac{\hat{Y}_{t,T}}{Y_t}<x\;\Big|\; \beta_t=j\right)\leq 1, 
 \quad (t,x,j)\in[0,T]\times[1,\infty)\times{\cal M}, 
\end{equation} 
 cf. the definition \eqref{g2} of $G$.
 Next, by \eqref{g2} and Lemma~\ref{equationg}, 
 for any $(t,x,j)\in [0,T]\times [1,\infty)\times {\cal M}$, we find 
\begin{align} 
\nonumber 
 & \mathbb{L}G(t,x,j) = x\sigma^2(j)P\left(\frac{\hat{Y}_{t,T}}{Y_t}<x\;\Big|\;\beta_t=j\right)-\mu(j)E\left[\max \left (x, \hat{Y}_{t,T} / Y_t \right)\;\Big|\;\beta_t=j\right]
\\
\label{h2terms}
&=E\left[ 
 x\sigma^2(j) 
 {\bf1}_{ \left\{ \hat{Y}_{t,T} / Y_t <x\right\} } 
 -\mu(j)\left(x\vee\frac{\hat{Y}_{t,T}}{Y_t}\right)
 \;\Big|\;\beta_t=j\right]
\\
\nonumber 
&=E\left[x(\sigma^2(j)-\mu(j)){\bf1}_{\left\{ \hat{Y}_{t,T} / Y_t <x\right\}}\;\Big|\;\beta_t=j\right]
-E\left[\mu(j)\left(x\vee\frac{\hat{Y}_{t,T}}{Y_t}\right){\bf1}_{\left\{ 
 \hat{Y}_{t,T} / Y_t \geq x\right\}}
\;\Big|\;\beta_t=j\right],
\end{align} 
 which shows \eqref{djkladas}, and implies by \eqref{h2terms} that 
 $\mathbb{L}G(t,x,j)$ is nondecreasing and continuous 
 in $t\in [0,T]$ when $\mu(j)\geq 0$. 
\end{Proof} 
 The proof of the next lemma, which will be used 
 in Proposition~\ref{2spcase} below, 
 extends the argument of \cite{dutoit} page~993 to the
 regime-switching setting. 
\begin{lemma} 
\label{h<0} 
 We have 
$$\{(t,x,j)\in[0,T)\times[1,\infty) \times{\cal M} \ : \ \mathbb{L}G(t,x,j) < 0\}\subset C,$$
where $C=D^c$ is the continuation set.
\end{lemma}
\begin{Proof}
 By Lemma~\ref{mgenerator} below and Lemma~1 in \cite{yaozhangzhou} 
 we have 
\begin{equation} 
\label{app1}
E[G(s,X_s^{t,x},\beta_s)\;|\;\beta_t=j] 
=G(t,x,j)+E\left[\int_t^s 
 \mathbb{L}G(r,X_r^{t,x},\beta_r)dr \;\Big|\;\beta_t=j\right], 
\end{equation} 
 $s \in [t,T]$. 
 Assume now that 
 $(t,x,j)\in[0,T)\times[1,\infty)\times{\cal M}$ is such 
 that $\mathbb{L}G(t,x,j) <0$. 
 By the continuity of $\mathbb{L}G(t,x,j) $ with respect to $t$, there exists an open neighbourhood 
 $U\subset [0,T) \times[1,\infty)$ of $(t,x)$ 
 such that $\mathbb{L}G(s,y,j) <0$ for all $(s,y) \in U$. 
 Replacing $s$ in \eqref{app1} with the first exist time 
 $\tau_U$ of $U$ when $(X_s^{t,x},\beta_s)_{s\in[t,T]}$ 
 is started at $(x,j)$ at time $t$, 
 Relation~\eqref{app1} above shows by optional sampling that 
\begin{eqnarray*} 
E[G(\tau_U,X_{\tau_U}^{t,x},\beta_{\tau_U})\;|\;\beta_t=j]=G(t,x,j)+E\left[\int_t^{\tau_U}\mathbb{L}G(r,X_r^{t,x},\beta_r)dr\;\bigg|\;\beta_t=j\right].
\end{eqnarray*}
 Since $\tau_U>t$ a.s. and 
 $\mathbb{L}G(r,X_r^{t,x},\beta_r)<0$ when $r\in( t ,\tau_U)$, the 
 right hand side is strictly smaller than $G(t,x,j)$, while we have 
$$ 
 E[G(\tau_U,X_{\tau_U}^{t,x},\beta_{\tau_U})\;|\;\beta_t=j]\geq V(t,x,j), 
$$ 
 showing that $V(t,x,j)<G(t,x,j)$, which implies that $(t,x,j)\in C$.
\end{Proof}
 Next we derive the following continuity result 
 wich has been used in the proof of Proposition~\ref{3.1}.  
\begin{lemma}
\label{contlemma}
 For any $j \in {\cal M}$, the mapping
 $(t,x)\longmapsto V(t,x,j)$ is jointly continuous on $[0,T]\times [1,\infty)$.
\end{lemma}
\begin{Proof} We proceed in two steps.
 $(i)$ We show that the mapping
 $t\longmapsto V(t,x,j)$ is continuous on $[0,T]$
 for every fixed $x\geq 1$ and any $j \in {\cal M}$. 
 By \eqref{vtjx} we have 
\begin{eqnarray}
 \nonumber
 V(t,x,j)&=&\inf\limits_{t\leq\tau\leq T}E\left[\frac{(xY_t)\vee\hat{Y}_{t,T}}{Y_\tau}\,\bigg|\,\beta_t=j\right]\\
\nonumber 
&=&\inf\limits_{0\leq\tau\leq T-t}E\left[\frac{x\vee \e^{\max\limits_{0\leq r\leq T-t}((\mu(\beta_r)-
      \sigma^2(\beta_r) /2 )r
      +\sigma(\beta_r)B_r)}}
  {\e^{(\mu(\beta_\tau)-
      \sigma^2(\beta_\tau) /2 )\tau+\sigma(\beta_\tau)B_\tau}}\,\bigg|\,\beta_0=j\right]
 \\
\nonumber 
  & = & \inf\limits_{0\leq\tau\leq T-t}E\left[U(t,\tau)\,\bigg|\,\beta_0=j\right],
  \quad
  t\in[0,T], \ j\in{\cal M}, \ x\in[1,\infty), 
\end{eqnarray} 
 where 
 $$U(t,s):=\frac{x\vee \e^{\max\limits_{0\leq r\leq T-t}((\mu(\beta_r)-
     \sigma^2(\beta_r) / 2 )r+\sigma(\beta_r)B_r)}}
 {\e^{(\mu(\beta_s)-
     \sigma^2(\beta_s) / 2 )\tau+\sigma(\beta_s)B_s}},
\qquad
s,t\in[0,T]. 
$$
 For any $({\cal G}_s^t)_{s\in [t,T]}$-stopping time $\tau\in[0,T-t]$ we have
 \begin{align}
   \nonumber 
   0&\leq E\left[U(t,\tau)-U(t+s,\tau)\,\bigg|\,\beta_0=j\right]
   \\
\nonumber
&\leq \sqrt{
  E\left[\left(\e^{\max\limits_{0\leq r\leq T-t}((\mu(\beta_r)-
      \sigma^2(\beta_r) / 2 )r+\sigma(\beta_r)B_r)}
    -\e^{\max\limits_{0\leq r\leq T-t-s}((\mu(\beta_r)-
      \sigma^2(\beta_r) / 2 )r+\sigma(\beta_r)B_r)}\right)^2
\,\bigg|\,\beta_0=j\right] } 
\\
\nonumber
&\quad \times
\sqrt{
  E\left[\e^{-2(\mu(\beta_\tau)-
      \sigma^2(\beta_\tau) / 2 )\tau-2\sigma(\beta_\tau)B_\tau} \,\bigg|\,\beta_0=j\right] } 
\\
\nonumber
&\leq
\sqrt{
  E\left[\left(\e^{\max\limits_{0\leq r\leq T-t}((\mu(\beta_r)-
      \sigma^2(\beta_r) / 2 )r+\sigma(\beta_r)B_r)}
    -\e^{\max\limits_{0\leq r\leq T-t-s}((\mu(\beta_r)-
      \sigma^2(\beta_r) / 2 )r+\sigma(\beta_r)B_r)}\right)^2
\,\bigg|\,\beta_0=j\right] } 
\\
\nonumber
&\quad\times
\sqrt{
  \e^{ (T-t) \max\limits_{i\in{\cal M}}|3\sigma(i)-2\mu(i)| }E\left[\e^{-2\sigma^2(\beta_\tau)\tau-2\sigma(\beta_\tau)B_\tau} \,\bigg|\,\beta_0=j\right] } 
\\
\nonumber
&\leq \sqrt{
  E\left[\left(\e^{\max\limits_{0\leq r\leq T-t}((\mu(\beta_r)-
      \sigma^2(\beta_r) / 2 )r+\sigma(\beta_r)B_r)}
    -\e^{\max\limits_{0\leq r\leq T-t-s}((\mu(\beta_r)-
      \sigma^2(\beta_r) / 2 )r+\sigma(\beta_r)B_r)}\right)^2
\,\bigg|\,\beta_0=j\right] } 
\\
\label{ia}
&\quad\times \e^{ (T-t) \max\limits_{i\in{\cal M}}|3\sigma(i)-2\mu(i)| / 2 }, 
\end{align}
where we applied the optional sampling theorem.
Letting $s$ tend to $0$ on both sides of \eqref{ia}, we get 
\begin{equation*}
\lim\limits_{s\searrow 0}E\big[U(t+s,\tau) \mid \beta_0=j\big]=E\big[U(t,\tau) \mid \beta_0=j\big],
\end{equation*}
and since the convergence is uniform on
all $({\cal F}_s)_{s\in [0,T]}$-stopping times $\tau\in[0,T]$, we obtain 
\begin{eqnarray}
  \nonumber
  \lefteqn{
    \! \! \! \! \!  \!  \!
    \liminf\limits_{s\searrow 0}\inf\limits_{0\leq\tau\leq T-t-s}E\big[U(t+s,\tau)
      \mid \beta_0=j\big]
\geq 
\liminf\limits_{s\searrow 0}\inf\limits_{0\leq\tau\leq T-t}E\big[U(t+s,\tau)  \mid \beta_0=j\big]
  }
  \\
  \label{fjhksdf} 
  &= & \inf\limits_{0\leq\tau\leq T-t}\lim\limits_{s\searrow 0}E\big[U(t+s,\tau)  \mid \beta_0=j\big]
= \inf\limits_{0\leq\tau\leq T-t}E\big[U(t,\tau)  \mid \beta_0=j\big].
\end{eqnarray}
Next, according to Proposition \ref{3.1} there exists
an optimal $({\cal F}_s)_{s\in [0,T]}$-stopping time $\tau^*_t\in[0,T-t]$ such that
\begin{equation}\label{iistar}
\inf\limits_{0\leq\tau\leq T-t}E\big[U(t,\tau) \mid \beta_0=j\big]=E\big[U(t,\tau^*_t) \mid \beta_0=j\big],
\end{equation}
hence we have
\begin{eqnarray}
  \label{iia}
  \inf\limits_{0\leq\tau\leq T-t-s} E\big[U(t+s,\tau) \mid \beta_0=j\big] & \leq
  &
  \inf\limits_{0\leq\tau\leq T-t-s}E\big[U(t,\tau) \mid \beta_0=j\big]\\
\nonumber
   &\leq & E\big[U(t,\tau^*_t\wedge(T-t-s)) \mid \beta_0=j\big].
\end{eqnarray}
Since $U (t,s)$ is nonnegative for any $s,t\in[0,T]$, we have
\begin{eqnarray*}
  U(t,\tau^*_t\wedge(T-t-s)) & \leq & U(t,\tau^*_t)+U(t,T-t-s)
  \\
&= & U(t,\tau^*_t)+\frac{x\vee \e^{\max\limits_{0\leq r\leq T-t}((\mu(\beta_r)-
    \sigma^2(\beta_r) / 2 )r+\sigma(\beta_r)B_r)}}
{\e^{(\mu(\beta_{T-t-s})-
    \sigma^2(\beta_{T-t-s}) / 2 )(T-t-s)+\sigma(\beta_{T-t-s})B_{T-t-s}}}
\\
&\leq & U(t,\tau^*_t)+\frac{x\vee \e^{\max\limits_{0\leq r\leq T-t}((\mu(\beta_r)-
    \sigma^2(\beta_r) / 2 )r+\sigma(\beta_r)B_r)}}
  {\e^{\inf\limits_{i\in{\cal M},r\in[0,T-t]}(\mu(i)-
      \sigma^2(i) / 2 )r+\inf\limits_{i\in{\cal M},r\in[0,T-t]}(\sigma(i)B_r)}},
\end{eqnarray*}
which is integrable by \eqref{iistar}. 
 By the reverse Fatou Lemma we have 
 \begin{eqnarray}
   \nonumber
   \lefteqn{
     \!\!\!\!\!\!\!\!\!\!\!\!\!\!\!\!\!\!\!\!\!\!\!
          \!\!\!\!\!\!\!\!\!\!\!\!\!\!\!\!\!\!\!\!\!\!\!
          \!\!\!\!\!\!\!\!\!\!\!\!\!\!\!\!\!\!\!\!\!\!\!
          \limsup\limits_{s\searrow 0}E\big[
            U(t,\tau^*_t\wedge(T-t-s)) \mid \beta_0=j\big]
     \leq
     E\big[\limsup\limits_{s\searrow 0}U(t,\tau^*_t\wedge(T-t-s)) \mid \beta_0=j\big]
   }
   \\
   \label{iib}
&= & E\big[U(t,\tau^*_t) \mid \beta_0=j\big].
\end{eqnarray} 
 Combining \eqref{iistar}, \eqref{iia}, \eqref{iib}
 and \eqref{fjhksdf} we find
 \begin{equation*}
\lim\limits_{s\searrow 0}\inf\limits_{0\leq\tau\leq T-t-s}E\left[U(t+s,\tau)\,\bigg|\,\beta_0=j\right]
=\inf\limits_{0\leq\tau\leq T-t}E\left[U(t,\tau)\,\bigg|\,\beta_0=j\right].
\end{equation*}
 Similarly we have
\begin{equation*}
\lim\limits_{s\searrow 0}\inf\limits_{0\leq\tau\leq T-t+s}E\left[U(t-s,\tau)\,\bigg|\,\beta_0=j\right]
=\inf\limits_{0\leq\tau\leq T-t}E\left[U(t,\tau)\,\bigg|\,\beta_0=j\right],
\end{equation*}
hence $t\longmapsto V(t,x,j)$ is continuous on $[0,T]$.
\\
\vspace{-0.4cm}
 
\noindent $(ii)$
We show that $x\longmapsto V(t,x,j)$ is continuous on $[1,\infty)$, uniformly in $t \in[0,T]$, 
 extending the argument of \cite{dutoit} page~995 
 to the regime-switching setting. 
 By Relation~\eqref{G'x} and the mean value theorem,
 for all  $y\in[x,\infty)$ there exists a (random)
 $\eta \in [X_{t+\tau}^{t,x}, X_{t+\tau}^{t,y}]$ such that
 for any $({\cal F}_s)_{s\in [0,T]}$-stopping time $\tau\in[0,T-t]$
 we have 
 \begin{eqnarray}
   \nonumber
   \lefteqn{
     \!\!\!\!\!\!\!\!\!\!\!\!\!\!\!\!\!\!\!\!\!\!\!\!\!\!\!\!\!\!\!\!\!\!\!\!\!\!\!\!\!\!
     \!\!\!\!\!\!\!\!\!\!\!\!\!\!\!\!\!\!\!\!\!\!\!\!\!\!\!\!\!\!\!\!\!\!\!\!\!\!\!\!\!\!
         G(t+\tau,X_{t+\tau}^{t,y},\beta_{t+\tau})-G(t+\tau,X_{t+\tau}^{t,x},\beta_{t+\tau})
     = \frac{\partial G}{\partial x}(t+\tau, \eta)(X_{t+\tau}^{t,y}-X_{t+\tau}^{t,x},\beta_{t+\tau})
   }
   \\
   \label{vcon2}
 & \leq & (y-x) \frac{Y_t}{Y_{t+\tau}},
\end{eqnarray} 
 since $X_{t+\tau}^{t,y}-X_{t+\tau}^{t,x}\leq {(y-x)Y_t}/{Y_{t+\tau}}$ by
 \eqref{x1}.
 Let now $(t,x,j)\in[0,T]\times[1,\infty)\times{\cal M}$
 and consider $\tau_x:=\tau(t,x,j)$ given by \eqref{tau}.
 By Lemma~\ref{v=g} we have
 \begin{equation}
   \label{vcon1}
V(t,y,j)-V(t,x,j)\leq E\left[
 G(t+\tau_x,X_{t+\tau_x}^{t,y},\beta_{t+\tau_x})-G(t+\tau_x,X_{t+\tau_x}^{t,x},\beta_{t+\tau_x})\;|\;\beta_t=j \right].
\end{equation}
 Since $E\big[ Y_t / Y_{t+\tau} \mid \beta_t=j\big]$
 is uniformly bounded as in \eqref{<infty},
 taking expectation on both sides of \eqref{vcon2} yields
\begin{equation}\label{vcon3}
\lim\limits_{y\to x}E \left[
 G(t+\tau,X_{t+\tau}^{t,y},\beta_{t+\tau})-G(t+\tau,X_{t+\tau}^{t,x},\beta_{t+\tau})\;|\;\beta_t=j\right]=0,
\end{equation}
 uniformly in $t \in [0,T]$ and in the $({\cal F}_s)_{s\in [0,T]}$-stopping times
 $\tau \in [0,T-t]$.
 Since $V(t,x,j)$ is increasing in $x\in[1,\infty)$,
 \eqref{vcon1} and \eqref{vcon3} yield
$$0\leq \lim\limits_{y\to x} (V(t,y,j)-V(t,x,j))\leq 0,$$
 which shows the continuity of $x\longmapsto V(t,x,j)$,
 uniformly in $t\in [0,T]$, for all $j \in {\cal M}$.
 \\

 \noindent
 From $(i)$ and $(ii)$ we conclude to the joint continuity of $(t,x)\longmapsto V(t,x,j)$ on $[0,T]\times [1,\infty)$ by classical arguments.    
\end{Proof} 
\begin{lemma}
\label{contlemma1}
 The mapping
 $(t,x)\longmapsto G(t,x,j)$ is jointly continuous on $[0,T]\times [1,\infty)$.
\end{lemma}
\begin{Proof}
  By Relation~\eqref{G'x} and the mean value theorem,
 for all $y\in[x,\infty)$ there exists an
 $\eta \in [x,y]$ such that
 for any $t\in[0,T]$
 we have
 $$0\leq G(t,y,j)-G(t,x,j)=
 (y-x)
  \frac{\partial G}{\partial x}(t,\eta ,j)
 \leq y-x,$$
 which shows the continuity of $x\longmapsto G(t,x,j)$, uniformly in $t\in[0,T]$. On the other hand, we have by \eqref{g2} that $t\longmapsto G(t,x,j)$ is continuous on $[0,T]$ for every $x \geq 1$. 
 We conclude to the joint continuity of $(t,x)\longmapsto G(t,x,j)$ on $[0,T]\times [1,\infty)$ by a classical argument. 
\end{Proof} 
 We close this section with the following lemma. 
\begin{lemma}\label{mgenerator}
 The Markov process $(s,X_s^{t,x},\beta_s)_{s\in[t,T]}$ has the 
 infinitesimal generator
\begin{equation*}
\mathbb{L}f(s,y,j)=\left(\frac{\partial}{\partial s}+y(\sigma^2(j)-\mu(j))\frac{\partial}{\partial y}+\frac{1}{2}\sigma^2(j)y^2\frac{\partial^2 }{\partial y^2}\right)f(s,y,j)+\sum_{i=1}^mq_{j,i}f(s,y,i),
\end{equation*} 
 $s\in [0,T]$, $j\in {\cal M}$, $y \in [1,\infty)$, 
 for $f\in \Dom (\mathbb{L})$ satisfying $\frac{\partial f}{\partial y}(s,1+,j)=0$. 
\end{lemma}
\begin{Proof}
 Letting
\begin{equation*}
 Z_s^{t,x}:= 
\log x\vee\hat S^t_{s-t}-\int_0^{s-t}\sigma(\beta_{u+t})d\hat B^t_u,
\end{equation*}
 $s\in [t, T]$, $x\geq 1$, from \eqref{xtx} we have 
 $X_s^{t,x}= \exp \left( Z_s^{t,x} \right)$, 
 $s\in [t, T]$, $x\geq 1$.
 Since $(\hat S^t_r)_{r\in [0,T-t]}$ is nondecreasing it has finite variation,
 hence
$$ d\langle Z_r^{t,x},Z_r^{t,x} \rangle = \sigma^2(\beta_r)\langle
d \hat B^t_{r-t},d\hat B^t_{r-t}\rangle = \sigma^2(\beta_r)d \langle B_r,B_r \rangle =\sigma^2(\beta_r)dr,
$$
 which shows that
\begin{align}\label{dx}
dX_s^{t,x}
&=X_s^{t,x}dZ_s^{t,x}+\frac{1}{2}X_s^{t,x}d\langle Z_s^{t,x},Z_s^{t,x} \rangle
\nonumber\\
&=X_s^{t,x}dZ_s^{t,x}+\frac{1}{2}\sigma^2(\beta_s)X_s^{t,x}ds\nonumber\\
&=X_s^{t,x}d (\log x\vee\hat S^t_{s-t})-\sigma(\beta_s) X_s^{t,x}d\hat B_{s-t}^t+\frac{1}{2}\sigma^2(\beta_s)X_s^{t,x}ds.
\end{align} 
 Given that
 $\frac{\partial f}{\partial y}(s,1+,j)=0$ for $(s,y,j)\in[0,T]\times[1,\infty)\times{\cal M}$,
 we have
\begin{eqnarray*}
 \frac{\partial f}{\partial y} (s,X_s^{t,x},\beta_s)
 \;d (\log x\vee\hat S^t_{s-t})
 & = & \frac{\partial f}{\partial y}(s,X_s^{t,x},\beta_s){\bf1}_{\{X_s^{t,x}>1\}}d(\log x\vee\hat S^t_{s-t})
\\ 
 & = & \frac{\partial f}{\partial y}(s,X_s^{t,x},\beta_s){\bf1}_{\{Z_s^{t,x}>0\}}d(\log x\vee\hat S^t_{s-t}) 
\\ 
 & = & 0,
\end{eqnarray*}
 since $d(\log x\vee\hat S^t_{s-t})=0$
 when $Z_s^{t,x}>0$, $s\in[t,T]$.
 From \eqref{dx} this shows that
$$\frac{\partial f}{\partial y}(s,X_s^{t,x},\beta_s)\;dX_s^{t,x}=\frac{\partial f}{\partial y}(s,X_s^{t,x},\beta_s)\left(-\sigma(\beta_s) X_s^{t,x}d\hat B_{s-t}^t+\frac{1}{2}\sigma^2(\beta_s)X_s^{t,x}dr\right),$$
 and we conclude the proof by It\^o's calculus.
\end{Proof}
\section{Solution of the free boundary problem}
\label{ss4} 
\begin{figure}[ht!]
  \hskip-0.2cm
  \begin{subfigure}{.4\textwidth}
  \centering
  \includegraphics[width=7.8cm,height=4.5cm]{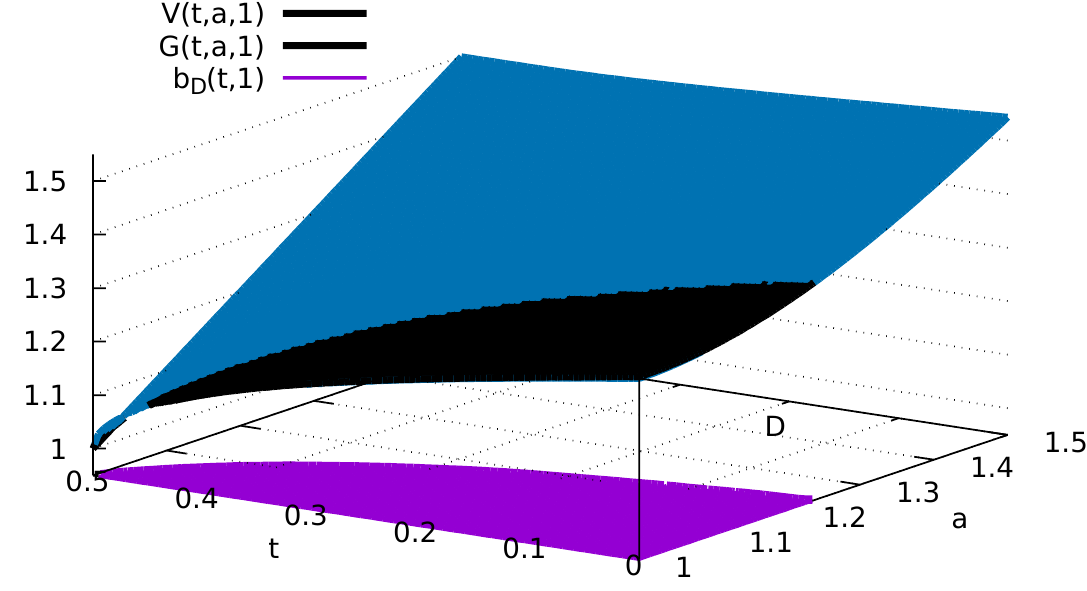}
\end{subfigure}
\hskip2cm
\begin{subfigure}{.4\textwidth}
  \centering
  \includegraphics[width=7.8cm,height=4.5cm]{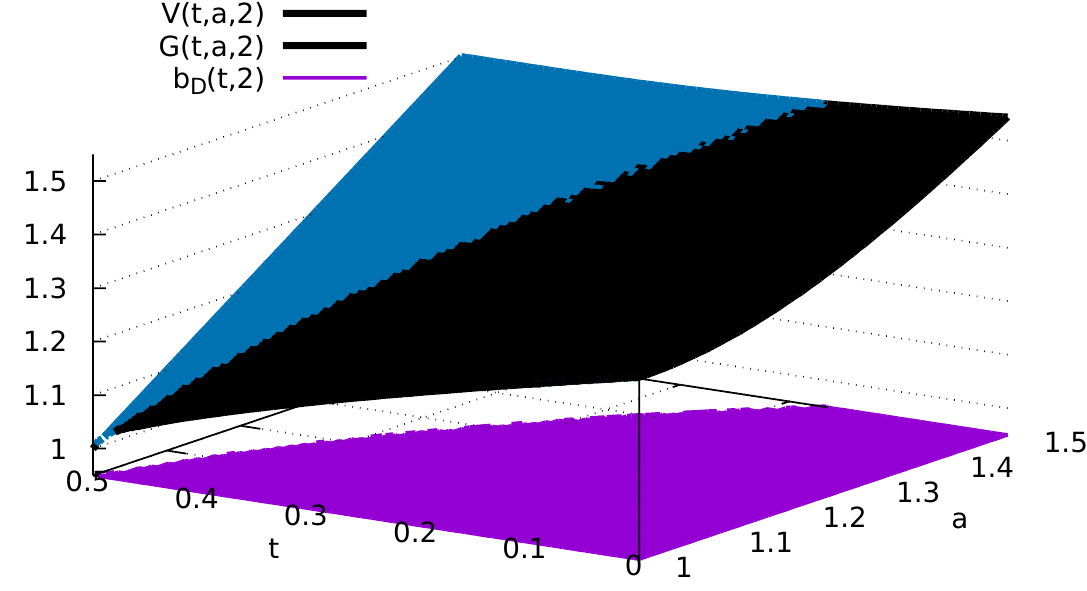}
\end{subfigure}
\vskip-0.2cm
\caption{Value functions in the two-state case.}
\label{fig2-3}
\end{figure}

 In this section we turn to the solution of the free boundary
 problem \eqref{vtjx}. 
 We start by providing sufficient conditions 
 on the drift coefficients $(\mu(j))_{j\in{\cal M}}$ 
 for the boundary function $b_D(t,j)$ defined by \eqref{bd} 
 to be nonincreasing and continuous in $t\in[0,T]$. 
 The next proposition~\ref{boundarydec} relies on Lemma~\ref{Ft} below.
\begin{prop} 
\label{boundarydec}
 Assume that $\mu(j)\geq 0$ for all $j\in{\cal M}$. 
 Then the boundary function $b_D(t,j)$ defined by \eqref{bd} is 
 nonincreasing in $t\in[0,T]$ and 
 continuous in $t\in[0,T]$, for all $j\in{\cal M}$. 
\end{prop}
\begin{Proof}
$(i)$ Monotonicity. 
 Let $(t,x,j) \in D$ and $s \in [t,T]$. 
 We have 
 $F(t,x,j)=0$ and $F(s,x,j)=0$ 
 since $F(t,x,j)$ is nondecreasing in $t$ by 
 Lemma~\ref{Ft}, hence 
$$
 [t,T] \times \{ x \} \times \{ j \} \subset D, 
$$
 showing that
 $(t,x,j)\in D\Longleftrightarrow [t,T] \times \{x\} \times \{ j \} \subset D$. 
 Then for any $s\in(t,T]$, we have $(s,b_D(t,j),j)\in D$ since $(t,b_D(t,j),j)\in D$. 
 By Proposition~\ref{3d} and noting that $(s,b_D(s,j),j)\in D$, we conclude 
 that $b_D(s,j)\leq b_D(t,j)$.
\\ 
\noindent 
$(ii)$ Right continuity. 
 Given $(t,b_D(t,j),j) \in D$, consider a strictly 
 decreasing sequence $(t_n)_{n \geq 1}$ such that 
 $\lim\limits_{n\to \infty}t_n=t$. 
 By part~$(i)$ above 
 we know that $b_D(t_n,j)\leq b_D(t,j)$, $n \geq 1$, and 
 $\lim\limits_{n\to \infty}b_D(t_n,j)\leq b_D(t,j)$.
 Next, by Proposition~\ref{3d} we have 
 $$
 [t,T]\times[b_D(t,j) ,\infty) \times\{j\} \subset 
 D,$$
 and since $(t_n,j,b_D(t_n,j))\in D$, $n \geq 1$, 
 and $D$ is closed, we have 
 $\left( 
 t, \lim\limits_{n\to \infty}b_D(t_n,j) ,j\right) \in D$, 
 hence $\lim\limits_{n\to \infty}b_D(t_n,j) \geq b_D(t,j)$. 
\\ 
\noindent 
$(iii)$ Left continuity. 
 Using Lemma~\ref{h<0} we can repeat the argument of \cite{dutoit} page~998 
 provided we show that the function $h(t,j)$ defined by 
\begin{equation} 
\label{htj}
h(t,j):=\inf\{x\in[1,\infty) 
 \ : \ 
 \mathbb{L}G(t,y,j)\geq 0, \quad 
 \forall y \in [x,\infty ) \}, 
\end{equation} 
 is continuous in $t\in [0,T]$ for all $j\in {\cal M}$,
 with $h(T,j) = 1$. 
 By Lemma~\ref{equationg} the function 
 $\mathbb{L}G(t,x,j)$ is nondecreasing in $t$ for all $x\geq 1$ 
 since $µ(j)\geq 0$ and it follows from the definition \eqref{htj} 
 of $h(t,j)$ that $t\longmapsto h(t,j)$ is nonincreasing in $t\in [0,T]$. 
 For any $t_0\in[0,T)$ and 
 decreasing sequence $( t_n )_{n\geq 1}\subset (t_0,T]$ 
 converging $t_0$ from the righthand side 
 we have 
 $\lim\limits_{n\to \infty}h(t_n,j)\leq h(t_0,j)$ 
 and 
 $\lim\limits_{n\to \infty}h(t_n,j)\geq h(t_k,j)$ for any $k\geq 1$, 
 hence $\lim\limits_{n\to \infty}h(t_n,j)\geq h(t_0,j)$ 
 as by the continuity of $t\to \mathbb{L}G(t,x,j)$ we have
$$ 
 \mathbb{L}G(t_0,\lim\limits_{n\to \infty}h(t_n,j),j)
 =\lim\limits_{k\to \infty}\mathbb{L}G(t_k,\lim\limits_{n\to \infty}h(t_n,j),j)\geq0,
$$
 and this proves that $\lim\limits_{t\searrow t_0}h(t,j)= h(t_0,j)$. 
 On the other hand we have 
 $h(t_0-,j):=\lim\limits_{t\uparrow t_0}h(t,j)\geq h(t_0,j)$ for any $t_0\in[0,T]$, 
 $j\in{\cal M}$. 
 In case $h(t_0-,j)>h(t_0,j)$ we have $\mathbb{L}G(t_0,x,j)\geq 0$ 
 for all $x\in[h(t_0,j),\infty]$. 
 In addition, for any  
 $t\in[0,t_0)$ and $x\in[h(t_0,j),h(t_0-,j))$ we have 
 $\mathbb{L}G(t,x,j)<0$ since $h(t,j)\geq h(t_0-,j)$, 
 hence $\mathbb{L}G(t_0,x,j)=0$ for all $x\in[h(t_0,j),h(t_0-,j))$ by 
 the continuity of $t \longmapsto \mathbb{L}G(t,x,j)$. 
 By Lemma~\ref{equationg} we would have
$$x\sigma^2(j)\frac{\partial G}{\partial x}(t_0,x,j)=\mu(j)G(t_0,x,j), 
 \qquad x\in[h(t_0,j),h(t_0-,j)), 
$$ 
 which shows that 
 $G(t_0,x,j)=C(t_0,j)x^{\mu(j) / \sigma^2(j) }$, where $C(t_0,j)$ 
 depends only on $t_0$ and $j\in {\cal M}$. 
 This is a contradiction since 
 $\frac{\partial }{\partial x}G(t_0,x,j) 
 =P\left( \hat{Y}_{t_0,T} / Y_{t_0} <x\;|\;\beta_{t_0}=j\right) 
 = C(t_0,j) \mu(j) x^{-1 + \mu(j) / \sigma^2(j)} / \sigma^2(j)$ 
 for $x\in[h(t_0,j),h(t_0-,j))$ 
 cannot hold when $\mu(j)<\sigma^2(j)$, 
 and more generally $\hat{Y}_{t_0,T} / Y_{t_0}$ cannot 
 have a power law, even locally. 
\end{Proof}
Figure~\ref{fig2-3} illustrates the result of
Proposition~\ref{boundarydec} by applying the recursive algorithm of 
       \cite{liu-privault2} in order to plot the
       value functions $V(t,a,j)$ and $G(t,a,j)$.
       In Figure~\ref{fig2-3} we take the positive drifts
       $\mu (1) =0.15$, $\mu (2) =0.05$, with $\sigma (1) =0.5$,
 $\sigma (2) =0.3$, $T=0.5$, $n=100$, $\delta_n = T/n = 0.05$, 
 and 
\begin{displaymath}
\mathbf{Q} =
\left[ \begin{array}{cc}
-2.5 & 2.5
\\
2 & -2 \nonumber
\end{array}\right].
\end{displaymath}
  
\noindent

Figure~\ref{fig2-3} also allows us to visualize
the stopping set $D$ and the continuation set 
\begin{equation*}
 C = \big\{
 (t,a,j)\in[0,T] \times [1,\infty) \times {\cal M} \ : \ V(t,a,j) < G(t,a,j)
   \big\}.
\end{equation*} 
The numerical instabilities observed are due to the necessity
to check the equality $V(t,a,j) = G(t,a,j)$ when
$V(t,a,j)$ and $G(t,a,j)$ are very close to each other.

\begin{figure}[ht!]
\centering
\includegraphics[height=0.3\textwidth,width=0.8\textwidth]{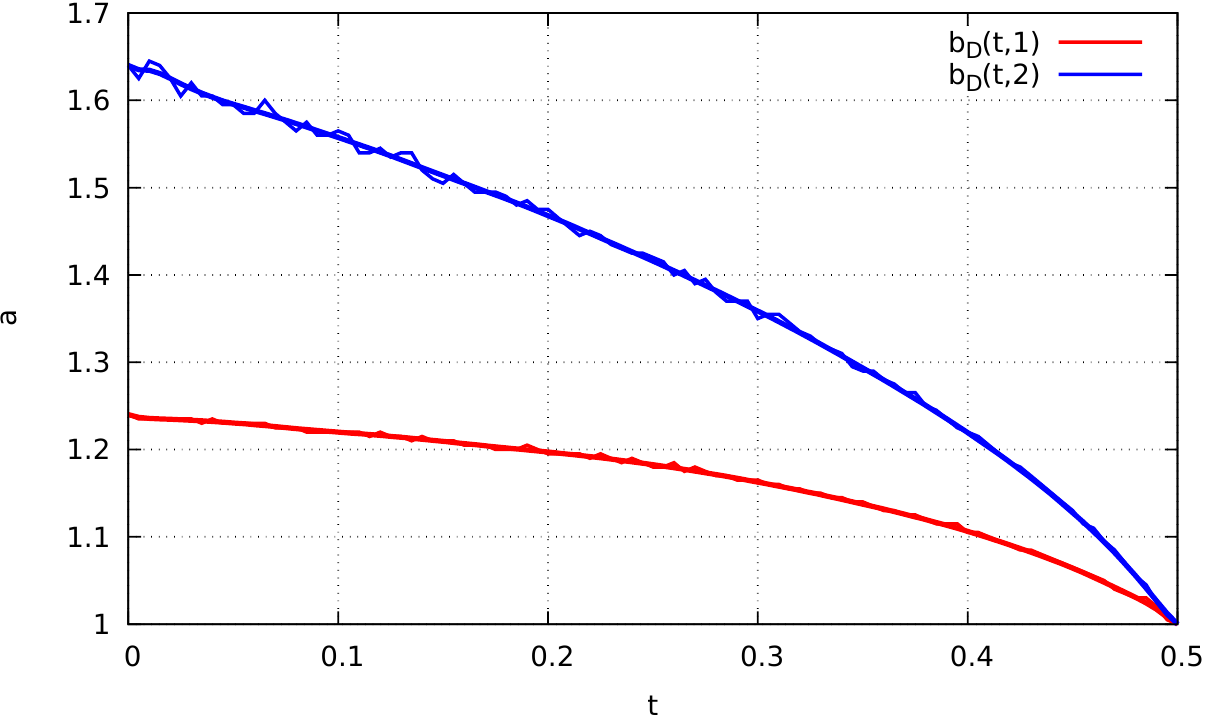}
\vskip-0.3cm
	\caption{Boundary functions in the two-state case.}
	\label{fig3}
\end{figure}

The boundary functions are plotted in Figure~\ref{fig3} 
based on Figure~\ref{fig2-3}, with spline smoothing. 
We observe that starting from state $1$ it is better to exercise
earlier than if we start from state $2$ which has a lower
drift.
This is due to the possibility to switch from
state $1$ to state $2$ after the average time 
$1/q_{1,1} = 0.4$ and to stay at state $2$ for the remaining
time $T-t \leq 1/q_{2,2} = 0.5$, 
in which case the drift takes the lower value $\mu (2)=0.05$.
The opposite occurs if we start from state $2$, for which the
boundary graph is higher than if we start from state $1$. 
\\
 
Similarly to (3.22)-(3.23) in \cite{dutoit}, 
 we now show that $F(t,x,j)$ defined by \eqref{discr} 
 is nondecreasing in $t\in [0,T]$ for all 
 $j\in {\cal M}$ and $x\in [1,\infty)$, 
   as in the following Lemma~\ref{Ft}
   which has been used for Proposition~\ref{boundarydec},
   and whose proof follows \cite{dutoit} 
 page~994. 
 \\

 We note that 
 without the condition $\mu(i)\geq 0$ for all $i\in{\cal M}$, 
 the function $F(t,x,j)$ in Lemma~\ref{Ft} may not be 
 nondecreasing in $t\in [0,T]$, 
 in which case the equivalence 
 $(t,x,j)\in D\Longleftrightarrow [t,T] \times \{x\} \times \{ j \} \subset D$ 
 in the proof of the next Proposition~\ref{boundarydec} 
 does not hold and in this situation 
 the boundary function $t\longmapsto b_D(t,j)$ may not 
 be decreasing in $t\in[0,T]$, cf. Figure~4 in \cite{liu-privault2}. 
\begin{lemma} 
\label{Ft}
 Under the condition 
 $\mu(j)\geq 0$ for all $j\in{\cal M}$, 
 the function 
$$ 
 t\longmapsto F(t,x,j)= V(t,x,j)-G(t,x,j)$$ 
 is nondecreasing in $t\in [0,T]$, 
 for any $(j,x)\in {\cal M}\times[1,\infty)$.
\end{lemma}
\begin{Proof} 
 For any $r,s\in[0,T-t]$, and $r<s$, denote $\tau_s:=\tau_D(s,x,j)-s\in[0,T-s]$ 
 by the definition \eqref{tau} of $\tau_D$.
 Replacing $s$ with $\tau_s$ and $t$ with $r$ in the formula \eqref{app1},
 and using optional sampling, we have
\begin{eqnarray} 
\label{F01}
F(r,x,j) & =&V(r,x,j)-G(r,x,j)
\\
\nonumber
& \leq&E[G(r+\tau_s,X_{r+\tau_s}^{r,x},\beta_{r+\tau_s})\;|\;\beta_r=j]-G(r,x,j)
\\
\nonumber
& =&E\left[\int_r^{r+\tau_s}\mathbb{L}G(v,X_v^{r,x},\beta_v)dv\;\Big|\;\beta_r=j\right]
\\
\nonumber
&=&E\left[\int_0^{\tau_s}\mathbb{L}G(v+r,X_{v+r}^{r,x},\beta_{v+r})dv\;\Big|\;\beta_r=j\right]
\\
\nonumber
& =&E\left[\int_0^{\tau_s}\mathbb{L}G(v+r,X_{v}^{0,x},\beta_{v})dv\;\Big|\;\beta_0=j\right].
\end{eqnarray}
 Combining \eqref{F01} with 
\begin{eqnarray*} 
  \lefteqn{
    \! \! \! \! \! \! \! \! \! \! \! \! \! \! \! \! \! \! \! \! \! \! \! \!
    \! \! \! \! \! \! \! \! \! \! \! \! \! \! \! \! \! \! \! \! \! \! \! \!
    F(s,x,j) = V(s,x,j)-G(s,x,j) 
    =
  E[G(s+\tau_s,X_{s+\tau_s}^{r,x},\beta_{s+\tau_s})\;|\;\beta_r=j]-G(s,x,j) 
 }
  \\
\nonumber
 & = & 
 E\left[\int_0^{\tau_s}\mathbb{L}G(v+s,X_{v}^{0,x},\beta_{v})dv\;\Big|\;\beta_0=j\right], 
\end{eqnarray*}
 we have
\begin{eqnarray} 
\label{F03}
\lefteqn{ 
 F(s,x,j)-F(r,x,j) 
} 
\\
\nonumber 
 & \geq & 
 E\left[\int_0^{\tau_s}\mathbb{L}G(v+s,X_{v}^{0,x},\beta_{v})dv\;\Big|\;\beta_0=j\right] 
 -E\left[\int_0^{\tau_s}\mathbb{L}G(v+r,X_{v}^{0,x},\beta_{v})dv\;\Big|\;\beta_0=j\right]
\\
\nonumber 
& =&E\left[\int_0^{\tau_s}\mathbb{L}G(v+s,X_{v}^{0,x},\beta_{v})-\mathbb{L}G(v+r,X_{v}^{0,x},\beta_{v})dv\;\Big|\;\beta_0=j\right].
\end{eqnarray} 
Since by \eqref{h2terms} the function
$t\mapsto \mathbb{L}G(t,x,i)$ is nondecreasing in $t$
when $\mu(i)\geq 0$, 
we find that the right hand side of \eqref{F03}
is nonnegative, thereby $F(t,x,j)$ is nondecreasing in $t\in[0,T]$.
\end{Proof}
\subsubsection*{Particular exercise strategies} 
 Next we show 
 that the stopping set $D$ has a simple form in 
 two special situations.
\begin{prop}
\label{2spcase}
 We have the following special cases of optimal stopping sets $D$.
\begin{enumerate}[i)]
\item Immediate exercise. Under the condition 
 $\mu(j)\leq 0$ for all $j\in{\cal M}$, 
 we have $D=[0,T]\times [1,\infty) \times{\cal M}$.
\item Exercise at maturity.
 Under the condition 
 $\mu(j)\geq \sigma^2(j)$ for all $j\in{\cal M}$,
 we have $D=\{T\}\times [1,\infty) \times{\cal M}$. 
\end{enumerate}
\end{prop}
\begin{Proof} 
 Replacing $s$ in \eqref{app1}
 with $\tau_D$ defined in 
 \eqref{tau} and using optional sampling, we find 
\begin{equation} 
\label{vge}
 V(t,x,i)=G(t,x,i)+E\left[\int_t^{\tau_{D(t,x,i)}} 
 \mathbb{L}G(r,X_r^{t,x},\beta_r)dr \;\Big| \; \beta_t = i \right], 
 \qquad 
 t\in [0,T].
\end{equation}
\begin{enumerate}[i)]
\item 
 In case $\mu(j)\leq 0$ for all $j\in{\cal M}$, by 
 Lemma~\ref{equationg}, 
 we have $\mathbb{L}G (t,x,i)>0$ for all $(t,x,i)\in [0,T)\times [1,\infty) \times{\cal M}$, hence \eqref{vge} implies $\tau_D(t,x,i)=0$ a.s., 
 otherwise it 
 contradicts the fact that $V(t,x,i) \leq G(t,x,i)$ 
 because of \eqref{vge}. 
 This implies $[0,T]\times [1,\infty) \times{\cal M} \subset D$. 
\item 
 In case $\mu(j)\geq \sigma^2(j)$ for all $j\in{\cal M}$, by 
 Lemma~\ref{equationg} 
 we have $\mathbb{L}G(t,x,i)<0$ for all $(t,x,i)\in [0,T)\times [1,\infty) \times{\cal M}$, and 
applying Lemma~\ref{h<0}, we see that $[0,T)\times [1,\infty) \times{\cal M}\subset C$, which means $D=\{T\}\times [1,\infty) \times{\cal M}$.
\vspace{-1cm} 
\end{enumerate}
\end{Proof}

 Finally we derive a Volterra type equation \eqref{ebdry} below 
 satisfied by the function $b_D(t,\beta_t)$
 defined in \eqref{bd}, for the boundary curves 
$$
 \left\{
 (t,x) \in [0,T] \times [1,\infty)
 \ : \ x = b_D (t,j)\right\}
$$
 of the optimal stopping set $D$ in
 \eqref{ddef}, for any $j \in {\cal M}$. 
\begin{prop}
\label{bdry}
 Assume that $\mu(j)\geq 0$, $j\in{\cal M}$. 
 The boundary function  $b_D (t,j)$ satisfies the
 Volterra
 type equation
\begin{equation} 
\label{ebdry}
G(t,b_D (t,j),j)=J(t,b_D (t,j),j)-\int_t^{T}K(t,r,b_D (t,j) ,j)dr,
\end{equation}
 $0 \leq t \leq T$,
 with terminal condition $b_D(T,j)=1$, 
 $j \in {\cal M}$,
 where
\begin{equation}\label{ebdry.2}
J(t,x,j):=E[X_{T}^{t,x}\;|\;\beta_t=j],
\end{equation}
 and
\begin{equation}\label{ebdry.3}
K(t,r,x,j):=E\left[\mathbb{L}V(r,X_r^{t,x},\beta_r){\bf 1}_{\{X_r^{t,x}> b_D(r,\beta_r)\}}\;\Big|\;\beta_t=j\right],
\end{equation}
 for $0\leq t\leq r\leq T$ and $x\geq 1$.
\end{prop}
\begin{Proof} 
 Noting that $V(t,x,j)\leq G(t,x,j)$ 
 for all $(t,x,j)\in[0,T]\times[1,\infty) \times{\cal M}$ by \eqref{vtj}, 
 the continuation set $C:= D^c$ is given by 
\begin{equation}\label{ccc}
C= D^c=\left\{(t,x,j)\in[0,T]\times [1,\infty) \times {\cal M} \ : \ V(t,x,j)<G(t,x,j)\right\}. 
\end{equation}
 According to Proposition~\ref{3.1}, for any $(t,x,j)\in C$, we have
\begin{equation}\label{vtau}
 V(t,x,j) =
 E \left[
 G \left( \tau_D, X_{\tau_D}^{t,x},\beta_{\tau_D}\right)\; \Big| \; \beta_t = j
 \right],
\end{equation}
 where $\tau_D=\tau_D(t,x,j)$ is defined by \eqref{tau}.
 Given that $\frac{\partial V}{\partial y} (t,1+,j)=0$ by Lemma~\ref{nr}, 
 by the application of \cite{peskir}, Chapter~III, \S~7.1.1, 
 \S~7.4.1 
 as in \cite{dutoit} \S~3.5, page~996,
 the function $V$ in \eqref{vtau} is ${\cal C}^{1,2}$ in the continuation set 
 $C$ in \eqref{ccc} 
 and it solves the Cauchy-Dirichlet free boundary problem
\begin{subequations}
\begin{empheq}[left=\empheqlbrace]{align}
\label{lv}
 & \displaystyle 
 \mathbb{L}V(t,y,j)= 0,\;\;\quad\quad\quad\quad \;(t,y,j)\in C,
\\ 
\nonumber 
\\ 
\label{lv2}
& V(t,y,j)=G(t,y,j),\quad\quad \;(t,y,j)\in \partial C,
\end{empheq}
\end{subequations}
 
 hence $\partial C\subset D$, 
 where $\partial C$ denotes the boundary of the open set $C$. 
 By the local time 
 change of variable formula of \cite{peskir2}, 
 and by Lemma~\ref{mgenerator} below with the property 
 $\frac{\partial V}{\partial y} (t,1+,j)=0$ shown in Lemma~\ref{nr} 
 above, we have 
\begin{eqnarray}
\nonumber
\lefteqn{
\lefteqn{E[X_T^{t,x}\;|\;\beta_t=j]=E[V(T,X_T^{t,x},\beta_T)\;|\;\beta_t=j]}\nonumber
}
\\
\label{p31a}
&=&V(t,x,j)+E\left[\int_t^T\mathbb{L}V(r,X_r^{t,x},\beta_r){\bf 1}_{\{X_r^{t,x}\neq b_D(r,\beta_r)\}}dr\;\Big|\;\beta_t=j\right]
\\
\nonumber
&& 
 \! \! \! \! \! \! \! \! \! \! \! \! 
 + \frac{1}{2} 
 E\left[\int_t^T\left( 
 \frac{\partial V}{\partial y} (r,X_r^{t,x}+,\beta_r) - 
 \frac{\partial V}{\partial y} (r,X_r^{t,x}-,\beta_r)\right)  {\bf 1}_{\{X_r^{t,x}=b_D(r,\beta_r)\}}d{\ell}^b_r(X^{t,x})\;\Big| \;\beta_t=j\right],
\end{eqnarray}
 where we applied the equality $V(T,X_T^{t,x},\beta_T)=X_T^{t,x}$, 
 and $({\ell}^b_r(X^{t,x}))_{r\in[t,T]}$ 
 denotes the local time of $X^{t,x}$ on the 
 (piecewise continuous and nonincreasing 
 by Proposition~\ref{boundarydec}) curve $r\longmapsto b_D(r,\beta_r )$. 
 By the smooth fit property shown in Lemma~\ref{sf} above, 
 the last term in \eqref{p31a} vanishes. 
 By Proposition~\ref{3d} above 
 and the definition \eqref{bd} of $b_D(t,j)$, 
 Relation~\eqref{lv} can be rewritten as
\begin{equation*} 
\mathbb{L}V(r,y,j ){\bf 1}_{\{ y <b_D(r,j)\}}=0, 
 \qquad r\in[0,T], \ j\in{\cal M}, \ y\geq 1, 
\end{equation*} 
 which implies 
\begin{equation}
\label{p31b}
E\left[\int_t^T\mathbb{L}V(r,X_r^{t,x},\beta_r)dr\;|\;\beta_t=j\right]=E\left[\int_t^T\mathbb{L}V(r,X_r^{t,x},\beta_r){\bf 1}_{\{X_r^{t,x}> b_D(r,\beta_r)\}}dr\;|\;\beta_t=j\right].
\end{equation} 
 Hence, combining \eqref{p31a} and \eqref{p31b}, we obtain
\begin{equation}\label{e=v+h}
E[X_T^{t,x}\;|\;\beta_t=j]=V(t,x,j)+\int_t^TE\left[\mathbb{L} V(r,X_r^{t,x},\beta_r){\bf 1}_{\{X_r^{t,x} > b_D(r,\beta_r) \}}\;|\;\beta_t=j\right]dr, 
\end{equation}
 and substituting $x$ with $b_D(t,j)$ in \eqref{e=v+h} above 
 we find that
\begin{align*}
\lefteqn{\!\!\!\!G(t,b_D (t,j),j)=V(t,b_D (t,j),j)}\\
&=E[X_T^{t,b_D (t,j)}\;|\;\beta_t=j]-E\left[\int_t^T\mathbb{L}V(r,X_r^{t,b_D (t,j)},\beta_r){\bf 1}_{\{X_r^{t,x}\geq b_D(r,\beta_r)\}}dr\;|\;\beta_t=j\right]\\
&=J(t,b_D (t,j),j)-\int_t^{T}K(t,r,b_D (t,j),j)dr,
\end{align*}
where the functions $J$, $K$ are defined by \eqref{ebdry.2}-\eqref{ebdry.3}.
\end{Proof}
\begin{remark}
  \label{r1}
  Note that the equation \eqref{ebdry} also involves 
the optimal value function $V(r,y,j)$
and not only the function $G(r,y,j)$. 
Indeed, when $m\geq 2$ 
the equality $V(r,y,j)=G(r,y,j)$ in \eqref{p31b}
for a given $(r,y,j) = (r,X_r^{t,x},\beta_r)\in D$ 
does not imply 
$$
\mathbb{L}V(r,y,j)=\mathbb{L}G(r,y,j) 
$$ 
as in \cite{dutoit} because we may not have $V(r,y,i)=G(r,y,i)$
for all $i=1,\ldots , m$ 
in the summation over the states of $(\beta_t)_{t\in[0,T]}$ in the 
definition \eqref{gene2} of $\mathbb{L}$.
In \cite{buffington} this issue is dealt with
via an ordering assumption 
on the boundary functions $(b_D(t,j))_{t\in [0,T]}$ 
in the two-state case $j=1,2$, see Assumption~3.1 therein,
however this method applies specifically to American options
and not to ultimate maximum problems, which have a more complex payoff
structure. 
Moreover, such an ordering condition may not be satisfied
in our current setting, cf. Figure~4 of \cite{liu-privault2}. 
\end{remark}
 In the absence of regime switching with 
 $Y_t=Y_0\e^{(\mu- \sigma^2 / 2 )t+\sigma B_t}$, 
 Relation~\eqref{bd} is replaced by 
 $$b_D (t) =\inf\{x \in \real_+  :(t,x)\in D\},\quad\quad t\in[0,T],
$$
 and the boundary equation \eqref{ebdry} becomes
\begin{equation*}
G(t,b_D(t))=E[X_{T}^{t,b_D(t)}]-E \left[\int_t^{T} 
 \mathbb{L} G( r,X_r^{b_D(t)}){\bf 1}_{\left\{X_r^{t,b_D(t)}> b_D(r)\right\}} dr \right],
\end{equation*}
 which recovers (3.50) in \cite{dutoit}, with 
$$ 
 \mathbb{L} G(r,x) 
 =\left(\frac{\partial}{\partial r}
 +x(\sigma^2 - \mu)\frac{\partial}{\partial x}+\frac{1}{2}\sigma^2 x^2\frac{\partial^2}{\partial x^2}\right)G(r,x),\quad\; r\in[0,T],\;x\in\real_+.
$$
 Since the Volterra type equation \eqref{ebdry} cannot
 be solved by standard methods under regime switching,
 we have applied the recursive algorithm of \cite{liu-privault2} 
 in order to plot Figures~\ref{fig2-3} and \ref{fig3}.
 
\footnotesize 

\def\cprime{$'$} \def\polhk#1{\setbox0=\hbox{#1}{\ooalign{\hidewidth
  \lower1.5ex\hbox{`}\hidewidth\crcr\unhbox0}}}
  \def\polhk#1{\setbox0=\hbox{#1}{\ooalign{\hidewidth
  \lower1.5ex\hbox{`}\hidewidth\crcr\unhbox0}}} \def\cprime{$'$}


\begin{thebibliography}{10}

\bibitem{buffington}
J.~Buffington and R.J. Elliott.
\newblock American options with regime switching.
\newblock {\em Int. J. Theor. Appl. Finance}, 5(5):497--514, 2002.

\bibitem{dutoit}
J.~du~Toit and G.~Peskir.
\newblock Selling a stock at the ultimate maximum.
\newblock {\em Ann. Appl. Probab.}, 19(3):983--1014, 2009.

\bibitem{eloe}
P.~Eloe, R.~H. Liu, M.~Yatsuki, G.~Yin, and Q.~Zhang.
\newblock Optimal selling rules in a regime-switching exponential {G}aussian
  diffusion model.
\newblock {\em SIAM J. Appl. Math.}, 69(3):810--829, 2008.

\bibitem{guoxin}
X.~Guo.
\newblock An explicit solution to an optimal stopping problem with regime
  switching.
\newblock {\em J. Appl. Probab.}, 38(2):464--481, 2001.

\bibitem{guozhang}
X.~Guo and Q.~Zhang.
\newblock Optimal selling rules in a regime switching model.
\newblock {\em IEEE Trans. Automat. Control}, 50(9):1450--1455, 2005.

\bibitem{hamilton}
J.D. Hamilton.
\newblock A new approach to the economic analysis of non-stationary time
  series.
\newblock {\em Econometrica}, 57:357--384, 1989.

\bibitem{le-wang}
H.~Le and C.~Wang.
\newblock A finite time horizon optimal stopping problem with regime switching.
\newblock {\em SIAM J. Control Optim.}, 48(8):5193--5213, 2010.

\bibitem{liu-zhang-yin}
R~H. Liu, Q.~Zhang, and G.~Yin.
\newblock Option pricing in a regime-switching model using the fast {F}ourier
  transform.
\newblock {\em J. Appl. Math. Stoch. Anal.}, pages Art. ID 18109, 22, 2006.

\bibitem{liu-privault2}
Y.~Liu and N.~Privault.
\newblock A recursive algorithm for selling at the ultimate maximum in
  regime-switching models.
\newblock Preprint, 18 pages, 2015.

\bibitem{peskir2}
G.~Peskir.
\newblock A change-of-variable formula with local time on curves.
\newblock {\em J. Theoret. Probab.}, 18(3):499--535, 2005.

\bibitem{peskir}
G.~Peskir and A.~Shiryaev.
\newblock {\em Optimal stopping and free-boundary problems}.
\newblock Lectures in Mathematics ETH Z\"urich. Birkh\"auser Verlag, Basel,
  2006.

\bibitem{Shiryaevstopping}
A.N. Shiryaev.
\newblock Optimal stopping rules.
\newblock {\em Springer-Verlag, New York, NY}, 1978.

\bibitem{shiryaev-xu-zhou}
A.N. Shiryaev, Z.~Xu, and X.Y. Zhou.
\newblock Thou shalt buy and hold.
\newblock {\em Quant. Finance}, 8(8):765--776, 2008.

\bibitem{vathana}
V.~Ly Vath and H.~Pham.
\newblock Explicit solution to an optimal switching problem in the two-regime
  case.
\newblock {\em SIAM J. Control Optim.}, 46(2):395--426, 2007.

\bibitem{yaozhangzhou}
D.D. Yao, Q.~Zhang, and X.Y. Zhou.
\newblock A regime-switching model for {E}uropean options.
\newblock {\em International Series in Operation Research and Management
  Science}, 94:281--300, 2006.

\end{thebibliography}
\end{document}